 \documentclass[12pt,reqno]{article}

 \usepackage{amsmath, amsthm, amsfonts}
 \usepackage{mathrsfs, amssymb}
 \usepackage{graphicx, color, bm, verbatim}
 \usepackage{pifont}

 \numberwithin{equation}{section}

 \newtheorem{thm}{Theorem}[section]
 \newtheorem{cor}[thm]{Corollary}
 \newtheorem{lem}[thm]{Lemma}
 \newtheorem{prop}[thm]{Proposition}
 \theoremstyle{remark}
 
 \theoremstyle{example}
 
\newtheorem{cond}[thm]{Condition}

 \def\ar{\!\!\!&}

 \def\mcr{\mathscr}\def\mbb{\mathbb}\def\mbf{\mathbf}
 \def\mrm{\mathrm}

 \def\beqlb{\begin{eqnarray}}\def\eeqlb{\end{eqnarray}}
 \def\beqnn{\begin{eqnarray*}}\def\eeqnn{\end{eqnarray*}}
 \def\d{{\mbox{\rm d}}}\def\e{{\mbox{\rm e}}}
 \def\blem{\begin{lem}\sl{}
 \def\elem{\end{lem}}}
 \def\bth{\begin{thm}\sl{}
 \def\eth{\end{thm}}}
 
 \def\bcond{\begin{cond}}\sl{}
 \def\econd{\end{cond}}
  \def\eqref#1{{\rm(\ref{#1})}}

 \def\black{\color{black}}

 \def\pf{\noindent{\it Proof.~~}}

\begin{document}

\noindent{(Version: \today)}
\
\bigskip\bigskip

\centerline{\Large\textbf{Exponential Ergodicity of CBIRE-Processes }}

\smallskip

\centerline{\Large\textbf{with Competition and Catastrophes}}

\bigskip

\centerline{Shukai Chen, \quad Rongjuan Fang, \quad Lina Ji, \quad Jian Wang}

\medskip

\bigskip

{\narrower{\narrower

\centerline{\textbf{Abstract}} We establish the exponential ergodic property in a weighted total variation distance of continuous-state branching processes with immigration in random environments with competition and catastrophes, under a Lyapunov-type condition and other mild assumptions. The proof is based on a Markov coupling process  along with some delicate estimates for the associated coupling generator. In particular, the main result indicates whether and how the competition mechanism, the environment and the catastrophe could balance the branching mechanism respectively to guarantee the exponential ergodicity of the process.

\medskip

\medskip

\par}\par}

\section{Introduction}
\subsection{Background}
{\it Continuous-state branching processes} (CB-processes) and {\it continuous-state branching processes with immigration} (CBI-processes) constitute important classes of Markov processes taking values in the positive half line. They were introduced as probabilistic models describing the evolution of large populations with small individuals. A stochastic equation for CBI-processes with branching mechanism $\Phi$ defined by
\beqlb\label{branching mechanism}
\Phi(\lambda)=b\lambda+\frac{1}{2}\sigma^2\lambda^2+\int_0^\infty(\e^{-\lambda z}-1+\lambda z)\,\mu(\d z),\quad\lambda\geq0
\eeqlb
was first established in \cite{DL06}. More explicitly,  a CBI-processes process is the unique strong solution to the following stochastic equation
\begin{equation}\label{SDE CBI}
X_t = X_0 -b\int_0^t X_s \,\d s + \sigma\int_0^t \sqrt{X_s} \,\d W_s + \int_0^t\int_0^\infty\int_0^{X_{s-}}z \,\tilde{M}(\d s, \d z, \d u)+I_t,
\end{equation}
where $b\in\mbb{R}$, $\sigma\ge0$, $(W_t)_{t\geq0}$ is a standard Brownian motion, $M(\d s, \d z, \d u)$ is a  Poisson random measure on $(0, \infty)^3$ with intensity $\d s\mu(\d z)\,\d u$ satisfying $\int_0^\infty(z\wedge z^2)\mu(\d z)<\infty$ and $\tilde{M}(\d s, \d z, \d u) = M(\d s, \d z, \d u) - \d s\,\mu(\d z)\,\d u$, and $(I_t)_{t\geq0}$ is a subordinator. A stochastic flow of discontinuous CB-processes was constructed in \cite{BL06} by using weak solutions of a special case of
(\ref{SDE CBI}).  We refer to \cite{DL12,FL10,LLWZ23,LM08} for more results of the stochastic equations of CBI-processes and stochastic flows of CBI-processes.

In recent years, the study of ergodicity for CB-processes and their generalized models has attracted considerable interest. The well used tools are the coupling approach and the  Meyn-Tweedie approach. For stochastic equations of the form (\ref{SDE CBI}), \cite{LM15} studied the exponential convergence in the total variation distance under the so-called Grey's condition by a coupling approach, see also \cite{FJKR23} in the Wasserstein distance. A natural generalization of the CBI-process is the so-called {\it affine Markov process}, which has also been used a lot in mathematical finance; see, e.g., \cite{DFS03}. In the affine framework, the Meyn-Tweedie approach applies to study the exponential ergodicity, see \cite{JKR17,MSV20,ZG18}. On the other hand, by applying the coupling methods, \cite{BW23,FJR20} and \cite{CL23} considered the exponential ergodicity in the Wasserstein distance and in the total variation distance, respectively.

Another generalized model is the interacting branching process to describe competitions or cooperations
among each pair of individuals in the population. \cite{LYZ19} considered a stochastic equation with the following form
\begin{equation}\label{sde nonlinear}
X_t=X_0+\int_0^t\gamma_0(X_s)\,\d s+\int_0^t\sqrt{\gamma_1(X_s)}\,\d W_s
+\int_0^t\int_0^\infty\int_0^{\gamma_2(X_{s-})}z\,\tilde{M}(\d s,\d z,\d u),
\end{equation}
where $\gamma_i\ (i=0,1,2)$ are continuous functions on $\mbb{R}_+$ satisfying certain assumptions, see \cite[Section 3]{LYZ19}. Clearly, (\ref{sde nonlinear}) includes CB-processes \cite{DL06}, logistic branching processes \cite{La05}, CB-processes with competition \cite{pardoux16} and so on. By making full use of the Markov coupling technique, \cite{LW20} obtained the exponential ergodicity in both the $L^1$-Wasserstein distance and the total variation norm, where the drift term is dissipative for large distance.  Under coupling methods and a Lyapunov-type condition inspired by \cite{DMT95,MT92,MT93a,MT93b}, \cite{LLWZ23} further studied the exponential ergodicity in a weighted total variation distance in the full range of criticality.

Branching processes in random environments were first introduced and studied in \cite{smithwilkinson}, where individuals in different generations may have different reproduction distributions. Those processes are more realistic compared with classical ones. From the mathematical point of view, they possess many interesting properties, such as the phase transition in the subcritical regime. We refer to \cite{BS15,LLGW} and references therein for more discussions. {\it Continuous-state branching processes in L\'evy random environments} (CBRE-processes) were introduced by
\cite{HLX18}, see also \cite{PP17}. Under certain moment condition on L\'evy environments,
such process is a strong solution of the following stochastic equation
\begin{equation}
X_t = X_0 -b\int_0^t X_s \d s + \sigma\int_0^t \sqrt{X_s} \d W_s + \int_0^t\int_0^\infty\int_0^{X_{s-}}z \tilde{M}(\d s, \d z, \d u)+ \int_0^t X_{s-} \d L_s,
\end{equation}
where $(L_t)_{t\geq0}$ is a L\'evy process defined by
$$
L_t=\beta_0t+\beta_1B_t+\int_0^t\int_{\mbb{R}}(\e^z-1)\,\tilde{N}(\d s,\d z),\quad t>0.
$$
Here $\beta_0\in \mbb{R}$, $\beta_1\ge0$, $(B_t)_{t\geq0}$ is a standard Brownian motion, and $\tilde{N}$ is a compensated Poisson measure  with intensity $\d s\nu(\d z)$ satisfying $\int_{-1}^1 z^2\nu(\d z) + \int_{|z| > 1} |\e^z - 1| \nu(\d z) < \infty.$  \cite{FJKR23} established the exponential ergodicity both in the $L^1$-Wasserstein distance and in the total variation norm with aid of a coupling perspective.

To model the evolution of the cell division with parasite infection, \cite{BT11} introduced a branching dynamic system with a tree structure, where the quantity of parasites in a cell evolves as a {\it Feller branching diffusion} (see \cite{MS20} for general CB-processes), the cell divides into two daughters in continuous time at a rate which may depend on the quantity of parasites, and parasites in the cells will be distributed in a random fraction. \cite{BT11} gave a criteria to determine whether the proportion of infected cells recovers or not. Note that the evolution of the quantity of parasites in a cell line plays a crucial role; see \cite[(3.1)]{BT11}. Inspired of this, \cite{MS21} introduced a generalized nonlinear CB-process, as a strong solution to
\beqlb\label{sde catastrophes}
X_t\ar=\ar X_0+\int_0^t\gamma_0(X_s)\,\d s+\int_0^t\sqrt{\gamma_1(X_s)}\,\d W_s
+\int_0^t\int_0^\infty\int_0^{\gamma_2(X_{s-})}z\,\tilde{M}(\d s,\d z,\d u)\cr
\ar\ar +\int_0^t\int_0^1\int_0^{r(X_{s-})}(z-1)X_{s-}\,Q(\d s,\d z,\d u),
\eeqlb
where $r(\cdot)$ is some nonnegative function on $\mbb{R}_+$, $Q$ is a Poisson random measure with intensity $\d sq(\d z)\d u$ with $q(\cdot)$ being a probability measure on $[0,1]$. Clearly, such process includes negative jumps and can be interpreted as the dynamics of the quantity of parasites in a cell line. Based on this, the last term on the right hand of (\ref{sde catastrophes}) is usually called the {\it catastrophic} part. It also can be seen as the state-dependent environment with pure negative jumps.

To the best of our knowledge, there are few known results on the exponential ergodicity in the total variation distance of generalized branching processes with negative jumps. The purpose of this paper is to establish the exponential ergodicity in a weighted total variation distance of {\it CBI-processes in random environments with competition and catastrophes} (CBIRE-processes with competition and catastrophes). For simplicity, we only consider a continuous immigration part determined by a drift coefficient $\alpha>0$. More precisely, with all the notations above at hand,  we assume that $(W_t)_{t\ge0}$, $(B_t)_{t\ge0}$, $\{M(\d s, \d z, \d u)\}$, $\{N(\d s, \d z)\}$ and $\{Q(\d s, \d z, \d u)\}$ are defined on a complete probability space and are independent of each other. Let us consider the following stochastic equation:
\beqlb\label{X}
X_t \ar=\ar X_0 + \int_0^t \gamma(X_s)\, \d s + \sigma\int_0^t \sqrt{X_s} \,\d W_s \cr
\ar\ar+ \int_0^t\int_0^\infty\int_0^{X_{s-}}z \tilde{M}(\d s, \d z, \d u)+ \int_0^t X_{s-} \d L^X_s,
\eeqlb
where $\gamma(x) = \alpha - bx- g(x)$ with $\alpha\ge0$ and $b\in \mbb R$, and
\begin{equation} \label{LX}
L^X_t=\beta_0t + \beta_1 B_t +\int_0^t\int_{\mbb{R}} (\e^z - 1) \tilde{N}(\d s, \d z)+ \int_0^t\int_0^1\int_0^{r(X_{s-})} (z - 1)Q(\d s, \d z,\d u)
\end{equation}with $\beta_0\in \mbb R$, $\beta_1\ge0$ and $r(x)$ being a nonnegative function on $\mbb R_+$.
Here $g$ is a {\it competition mechanism}, which by definition is a nondecreasing and continuous function on $\mbb{R}_+$ satisfying $g(0) = 0.$

\subsection{Main result}
To illustrate our main contributions, we present the following statement for the exponential ergodicity of the process $(X_t)_{t\geq0}$ determined by (\ref{X}). To do so, we first recall some necessary notation.

Given a Borel function $V\geq1$ on $\mbb{R}_+$, by ${\cal P}_V(\mbb{R}_+)$ we denote the space of all Borel probability measures $\varrho$ on $\mbb{R}_+$ satisfying
$$
\int_{\mbb{R}_+}V(x)\,\varrho(\d x)<\infty.
$$
Given $\pi_1,\pi_2\in{\cal P}_V(\mbb{R}_+)$, a coupling $H$ of $(\pi_1,\pi_2)$ is a Borel probability measure on $\mbb{R}_+\times\mbb{R}_+$ which has marginals $\pi_1$ and $\pi_2$, respectively. We write ${\cal H}_V(\pi_1,\pi_2)$ for the collection of all such couplings. Let $W_V$ be the $V$-{\it weighted total variation distance} between $\pi_1$ and $\pi_2\in{\cal P}_V(\mbb{R}_+)$ given by
 $$
W_V(\pi_1,\pi_2)=\int_{\mbb{R}_+}V(x)|\pi_1-\pi_2|(\d x),\quad \pi_1,\pi_2\in{\cal P}_V(\mbb{R}_+),
$$
where $|\cdot|$ denotes the total variation measure. We shall see that $W_V$ is actually the {\it Wasserstein distance} determined by the metric
\beqlb\label{defin dV}
d_V(x,y)=[V(x)+V(y)]
\mbf{1}_{\{x\neq y\}};
\eeqlb
that is,
$$
W_V(\pi_1,\pi_2)=\inf_{H\in{\cal H}_V(\pi_1,\pi_2)}\int_{\mbb{R}_+\times\mbb{R}_+}d_V(x,y)\,H(\d x,\d y).
$$
 We refer to \cite{C04} for the details. In particular, if $V\equiv1$, then $W_V$ reduces to the standard total variation distance.

Throughout this paper,
denote by $X:=(X_t)_{t\geq0}$ the unique strong solution of (\ref{X}). Let $P_t(x,\cdot)$ and $(P_t)_{t\geq0}$ be the transition function and the transition semigroup of the process $X$, respectively. We say the process $X$ is {\it exponentially ergodic in terms of the distance} $W_V$, if there are a unique stationary distribution $\pi$ and a constant $\lambda>0$ so that for all $\varrho\in{\cal P}_V(\mbb{R}_+)$ and $t>0$,
\beqlb\label{exponential ergodic}
W_{V}(\varrho P_t,\pi)\leq C(\varrho)\e^{-\lambda t},
\eeqlb where $\varrho\mapsto C(\varrho)$ is a nonnegative function  on ${\cal P}_V(\mbb{R}_+)$.

\bth\label{example}  Suppose that in the SDE \eqref{X} the rate function $r(x)$ is globally Lipschitz and $\alpha>0$. Let $(X_t)_{t\geq0}$ be a unique strong solution to (\ref{X}).  Let $V(x)=(x+1)^{\theta}$ with $\theta\in (0,1)$. Assume that
\beqlb\label{technial cond}
\limsup\limits_{x\rightarrow \infty}\frac{H(x)}{V(x)}=0,
\eeqlb
where
\beqnn
H(x):= \int_{0}^{1/x} (1 -  zx  )^3\big( \nu(\d \ln z)+  r(x)q(\d z)\big).
\eeqnn
Suppose that one of the following assumptions holds:

\begin{itemize}
\item[{\rm (i)}] $\sigma>0$;

\item[{\rm (ii)}] $\int_0^1z\,\mu(\d z)=\infty$ and there exist constants $c_0 > 0$ and $\delta > 0$ such that for all $|x|\le c_0$,
\beqnn
(\mu \wedge (\delta_x * \mu))(\mbb{R}_+)\ge \delta.
\eeqnn
\end{itemize}
Then the process $(X_t)_{t\geq0}$ is exponentially ergodic in terms of the distance $W_V$, if
\begin{equation}\label{condition of exam}
\begin{split}
	&\limsup_{x\rightarrow\infty}\Big[-\frac{g(x)}{x} + r(x)  \int_0^1(z^\theta - 1)q(\d z)\Big]\\
&~~~~~+{\beta_0}-b+{\frac{(\theta-1)\beta_1^2}{2}}+ \int_{-\infty}^\infty\left[\e^{z\theta} - 1 - \theta(\e^z - 1)\right]\nu(\d z) <0.
\end{split}
\end{equation}
\eth

 We now give some comments on Theorem \ref{example} and its proof.

 \begin{itemize}
 \item[{\rm(1)}] The condition that the constant $\alpha>0$ in the drift term $\gamma(x)$ roughly ensures that the stationary distribution of $(X_t)_{t\geq0}$ is not a degenerate distribution at zero, since, if $\alpha=0$, $X_t=0$ for all $t>0$ when $X_0=0$. Condition (i) or Condition (ii) in Theorem \ref{example} guarantees the existence of random perturbations of the branching part, which has been
      used in the study of exponential ergodicities in a weighted total variation norm of CBI-processes with competition, see \cite{LLWZ23}. Moreover, it is actually weaker than the corresponding assumptions for exponential ergodicities of Ornstein--Uhlenbeck type processes with nonlinear drift or nonlinear CB-processes in the total variation norm and in Wasserstein distances, see \cite{LW20,LW19} and references therein.
    Furthermore,  condition (\ref{condition of exam}) helps us to understand that whether and how the competition mechanism, the environment and the catastrophe could balance the branching mechanism to guarantee the exponential ergodicity. In particular, since $\e^{z\theta} - 1 - \theta(\e^z - 1)\le0$ for all $z\in \mbb{R}$ and $\theta\in (0,1)$, both the Gaussian noise and the Poisson noise in the catastrophe part  facilitate the exponential ergodicity of the process; on the other hand, from \eqref{condition of exam} we shall see that the drift term $\beta_0$ in the  environment part has the same status with the drift term of the branching part.

\item[{\rm(2)}]
The advantage of Theorem \ref{example} (or the general result Theorem \ref{main result} below) is that it not only works for general branching mechanisms without criticality restriction (see \cite{CL23,FJR20,JKR17,LM15,ZG18} for subcritical-type assumptions), but also is an effective exploration of the exponential ergodicity with explicit rates of branching processes with negative jumps (this is an essential difference between \cite{LLWZ23}, where only nonnegative jumps are considered, and the present paper). To handle the effect on the ergodicity arise from negative jumps,  we have to assume $\sigma>0$ or $\int_0^1z\,\mu(\d z)=\infty$ to guarantee
\beqlb\label{Phi_lambda1}
\Phi({ \lambda_1}) + \left(\inf_{x\geq0}(r(x))\int_0^1(1 - z)q (\d z) - \beta_0\right) { \lambda_1} > 0
\eeqlb
for some $\lambda_1>0$, which is crucial in the proof. Especially, if the environment and the catastrophe part vanish, (\ref{Phi_lambda1}) reduces to Condition (1.1) in \cite{LLWZ23}. In this case, (i) and (ii) can also reduce to Condition (1.2) in \cite{LLWZ23} (note that $\sigma>0$ or $\int_0^1z\,\mu(dz)=\infty$ is a necessary
condition for Grey's condition to hold).

\item[{\rm (3)}]The approach of Theorem \ref{example} is based on the refined basic coupling for the branching-jump term, the synchronous coupling for environment-jump term, the classic basic coupling for the catastrophic term and the coupling by reflection for Brownian motions. Namely, different couplings will apply different parts of the CBIRE-processes with competition and catastrophes, due to their different roles played in the ergodicity of the processes. Moreover, to efficiently realize the coupling idea to the CBIRE-processes with competition and catastrophes, we will use a suitable Lyapunov function to estimate the coupling generator on an unbounded area and we choose a nonsymmetric control function for the small distance.  For example, (\ref{technial cond}) is our technical condition, which is needed since the proofs of Theorem \ref{example} and Theorem \ref{main result} (see Section 3 for details) are based on a nonsymmetric control function associated with the weighted total variation distance. Though this approach partly is inspired by that in \cite{LLWZ23,LMW21}, we should carefully handle the negative jumps of our model. In particular, different from all the cited papers above, we can not use of the order-preservation property due to the fact that the CBIRE-processes with competition and catastrophes do not enjoy such kind property. To the best of our knowledge, this is the first time in the literature to study the ergodicity of branching processes (or their variants) with negative jumps.
\end{itemize}

The remainder of this paper is arranged as follows. In Section 2, we present some results on the existence and the uniqueness of the strong solution to (\ref{X}) and a Markov coupling of the unique strong solution through the construction of a coupling generator. General results on the exponential ergodicity of the strong solution to (\ref{X}) are stated in Section 3, where the proof of Theorem \ref{example} is given here.

\section{Unique strong solution and its coupling process}

This section consists of two parts. We first give the existence and the uniqueness of the strong solution to the SDE \eqref{X}, and then construct a new Markov coupling process for this unique strong solution.

\subsection{Existence and uniqueness of strong solution}

Let $C^2(\mbb{R}_+)$ be the linear space of twice continuously differentiable functions on $\mbb{R}_+.$ For any $f \in C^2(\mbb{R}_+)$,  write
\beqlb\label{L}
Lf(x) = L_bf(x) + L_ef(x)+L_cf(x),
\eeqlb
where
\beqnn
L_bf(x) \ar=\ar \gamma(x)f'(x) + \frac{1}{2}\sigma^2 xf''(x)+ x \int_0^\infty \Big[f(x + z) - f(x) - zf'(x)\Big]\,\mu(\d z),\cr
L_ef(x) \ar=\ar \beta_0xf'(x) + \frac{\beta_1^2x^2}{2}f''(x)+ \int_{-\infty}^\infty \Big[f(x\e^z) -f(x) -x(\e^z - 1)f'(x)\Big]\nu(\d z)
\eeqnn
and
\beqnn
L_cf(x) = r(x)\int_0^1\Big[f(z x) -f(x)\Big]q(\d z).
\eeqnn
Let $\mcr{D}(L)$ denote the linear space consisting of functions $f\in C^2(\mbb{R}_+)$ such that the three integrals involved in the operators $L_b$, $L_e$ and $L_c$ are convergent and define continuous functions on $\mbb{R}_+.$ In particular, $C_b^2(\mbb{R}_+) \subset \mcr{D}(L),$ where $C_b^2(\mbb{R}_+)$ is the space of bounded and continuous functions on $\mbb{R}_+$ with bounded and continuous derivatives up to the second order. However, in general for $f\in \mcr{D}(L)$, both $f$ and $Lf$ can be unbounded.
\bth\label{eus}
Suppose that $r(x)$ is locally Lipschitz on $\mbb{R}_+$. Then for any initial value $X_0 = x \ge 0$, there exists a unique nonnegative non-explosive
strong solution to the SDE \eqref{X}.
\eth

\pf The proof of the uniqueness of the strong solution to the SDE \eqref{X} is an application of \cite[Proposition 1]{PP18}, which is a combined result of \cite[Theorem 2.1]{DL12} and \cite[Theorem 6.1]{PL12}. Let $E=\{1,2\}$,
and
\beqnn
W(\d s,\d u)=\d W_s \delta_1(\d u) + \d B_s\delta_2(\d u),
\eeqnn
which is a white noise in $\mbb{R}_+\times E$ with intensity $\d s (\delta_1(\d u) +\delta_2(\d u))$. Let
\beqnn
\ar\ar U = [0,1)\times\mbb{R}_+,  \cr
\ar\ar V =\{1\}\times\mbb{R}_+ \times \mbb{R}_+ \cup \{2\}\times\mbb{R}_+\times\{0\},
\eeqnn
and
\beqnn
\ar\ar M_0(\d s,\d z,\d u):=Q(\d s,\d z,\d u),\cr
\ar\ar N_0(\d s,\d r,\d z,\d u):=\delta_1(\d r)M(\d s,\d z,\d u)+\delta_2(\d r)N(\d s,\d z)\delta_0(\d u),
\eeqnn
which are Poisson random measures.
Then \eqref{X} can be written as
\beqnn
X_t \ar=\ar X_0 + \int_0^t b(X_s)\d s + \int_0^t\int_E \sigma(X_s,u) W(\d s,\d u)+\int_0^t\int_U g_0(X_{s-},z,u) M_0(\d s,\d z,\d u)\cr
\ar\ar +\int_0^t\int_V h_0(X_{s-},r,z,u) \tilde{N}_0(\d s,\d r,\d z,\d u),
\eeqnn where
\beqnn
\ar\ar b(x) = \gamma(x)+\beta_0x,\cr
\ar\ar \sigma(x,u)=\sigma\sqrt{x}\mbf{1}_{\{u=1\}}+\beta_1 x\mbf{1}_{\{u=2\}},\cr
\ar\ar g_0(x,z,u)=(z-1)x\mbf{1}_{\{u<r(x)\}},\cr
\ar\ar h_0(x,r,z,u)=\mbf{1}_{\{r=1\}}z\mbf{1}_{\{u<x\}}+\mbf{1}_{\{r=2\}}x(\e^z-1)\mbf{1}_{\{u=0\}}.
\eeqnn
It is easy to verify that $(b,\sigma,g_0,h_0)$ are admissible and satisfy conditions (a), (b) and (c) in \cite[Proposition 1]{PP18}. Therefore, there exists a unique strong solution to the SDE \eqref{X}.

Now we prove that this unique strong solution $(X_t)_{t\ge0}$ is non-explosive. Let $\zeta_n := \inf\{t \ge 0: X_t \ge n\}$ for $n \ge 1.$ Note that $r(x)$ is nonnegative.
Similar to \cite[Proposition 2.3]{FL10},
we can obtain that there exists a constant $K\ge0$ such that
\beqnn
\mbb{E}[1 + X_{t\wedge \zeta_n}] \ar\le\ar \mbb{E}[1 + X_0] + \mbb{E}\left[\int_0^{t\wedge \zeta_n}(\gamma(X_s) + \beta_0X_s)\d s\right]\cr
\ar\le\ar \mbb{E}[1 + X_0] + K\mbb{E}\left[\int_0^{t\wedge \zeta_n}(1 + X_s)\d s\right].
\eeqnn
By Gronwall's lemma,
\beqnn
\mbb{E}[1 + X_{t\wedge \zeta_n}] \le \mbb{E}[1 + X_0]\exp\{Kt\}, \qquad t \ge 0.
\eeqnn
In particular,
$$
	(1 + n)\mbb{P}(\zeta_n \le t) \le \mbb{E}[1 + X_0]\exp\{Kt\}.
$$
holds since $X_{\zeta_n} \ge n$, and so one can see  that
\beqnn
\lim_{n \rightarrow \infty} \mbb{P}(\zeta_n \le t) = 0
\eeqnn
for any $t \ge 0$. It follows by Fatou's lemma that $\zeta_n \rightarrow \infty$ as $n \rightarrow \infty,$ which implies the result.
\qed

The next result justifies the fact that the operator $(L, \mcr{D}(L))$ defined by \eqref{L} is a {\it restriction of the generator} of the process $(X_t)_{t \ge 0}$. The proof is omitted here since it is similar to that of \cite[Proposition 4.2]{FL10}.

\begin{thm}
	Let $(X_t)_{t \ge 0}$ be the unique nonnegative non-explosive
	strong solution to \eqref{X}. Then for any $f \in \mcr{D}(L)$ and $n\ge1$,
	\beqnn
	f(X_{t\wedge\zeta_n}) - f(X_0) - \int_0^{t\wedge \zeta_n}Lf(X_s)\d s, \qquad t \ge 0
	\eeqnn
	 is a martingale, where   $\zeta_n := \inf\{t \ge 0: X_t \ge n\}$.
\end{thm}

\subsection{Markov coupling process}\label{subsection2.2}

In order to construct a Markov coupling of the process $(X_t)_{t \ge 0}$ determined by \eqref{X}, we begin with the construction of a new coupling operator for its generator $L$ given by \eqref{L}. Recall that $(X_t, Y_t)_{t \ge 0}$ is a Markov coupling of the process $(X_t)_{t \ge 0}$ given by \eqref{X}, if $(X_t, Y_t)_{t \ge 0}$ is a Markov process on $[0, \infty)^2$ such that the marginal process $(Y_t)_{t \ge 0}$ has the same transition probability as the process $(X_t)_{t \ge 0}.$ Denote by $\tilde{L}$ the infinitesimal generator of the Markov coupling process $(X_t, Y_t)_{t \ge 0}.$ Then the operator $\tilde{L}$ satisfies the following marginal property, i.e., for any $f_1, f_2 \in C^2(\mbb{R}_+),$
$$
\tilde{L}h(x, y) = Lf_1(x) + Lf_2(y),
$$
where $h(x, y) = f_1(x) + f_2(y)$ for any $x, y \in \mbb{R}_+,$ and $L$ is given by \eqref{L}. We call that $\tilde{L}$ is a {\it coupling operator} of $L.$

In this paper, as mentioned before, we will combine the refined basic coupling developed in \cite{LW19} for the branching-jump term, the synchronous coupling for environment-jump term, the classic basic coupling for the catastrophic term and the coupling by reflection for Brownian motions. Roughly speaking, the coupling by reflection for Brownian motion in the present setting means that we will take $(W_t)_{t\ge0}$ and $(-W_t)_{t\ge0}$ (resp., $(B_t)_{t\ge0}$ and $(-B_t)_{t \ge 0}$) for two marginal processes before they meet. For the jumping system driven by Poisson random measure $M$ we apply the refined basic coupling.
Then, the jumping system of the branching mechanism corresponding to the refined basic coupling of the operator $L$ is given by
\beqnn
(x, y) \rightarrow \left\{
\begin{aligned}
	&	(x+z, x + z), \qquad\qquad\qquad \frac12 (x\wedge y) \mu_{-(x - y)}(\d z),\\
	&	(x+z, 2y - x + z), \qquad\qquad \frac12 (x\wedge y)\mu_{(x - y)}(\d z),\\
	&	(x + z, y + z), \qquad\qquad\qquad  (x\wedge y)\Big[\mu(\d z) - \frac12 \mu_{-(x - y)}(\d z) - \frac12  \mu_{(x - y)}(\d z)\Big],\\
	&	(x + z, y),  \qquad\qquad\qquad\qquad (x - y)^+ \mu(\d z),\\
	&	(x, y+z),  \qquad\qquad\qquad\qquad (x - y)^- \mu(\d z),
\end{aligned}
\right.
\eeqnn
where $\mu_x(\d z) = (\mu \wedge (\delta_x * \mu))(\d z)$ for all $x \in \mbb{R}.$
We refer to \cite{LW19} for the details of refined basic couplings of L\'evy processes. Moreover, the jumping system driven by Poisson random measure $N$ of random environment corresponding to the synchronous coupling {\it w.r.t.} the jump size of the operator $L$ is given by
\beqnn
(x, y) \rightarrow (x\e^z, y\e^z), \quad\quad \nu(\d z);
\eeqnn
while, for the jumping system driven by $Q$ of catastrophe phenomenon, we use the basic coupling, i.e.,
\beqnn
(x, y) \rightarrow \left\{
\begin{aligned}
&	(z x, z y),  \qquad\qquad (r(x)\wedge r(y))q(\d z) ,\\
&  (z x, y),  \qquad\qquad~~(r(x)-r(y))^+q(\d z),
\\
&  (x, z y),  \qquad\qquad~~(r(x)-r(y))^-q(\d z).
\end{aligned}
\right.
\eeqnn
The readers can refer to \cite{C04} for the details of those two couplings.

Let $\Delta = \{(z, z): z \ge 0\} \subset \mbb{R}_+^2$ and $\Delta^c = \mbb{R}_+^2\backslash \Delta.$ With the aid of the idea above, given a function $f$ on $\mbb{R}_+^2$ twice continuously differentiable on $\Delta^c,$ we define
\beqlb\label{tildeL}
\tilde{L}f(x, y) = \tilde{L}_bf(x, y) + \tilde{L}_ef(x, y)+\tilde{L}_cf(x,y), \qquad (x, y) \in \Delta^c,
\eeqlb
where
\beqlb\label{tildeL0}
\tilde{L}_bf(x, y) \ar=\ar \gamma(x)f_x'(x, y) + \gamma(y)f_y'(x, y)\cr
\ar\ar + \frac{1}{2}\sigma^2 xf_{xx}''(x, y) + \frac{1}{2}\sigma^2 yf_{yy}''(x, y)  -\sigma^2\sqrt{xy}f_{xy}''(x, y)\cr
\ar\ar + \frac{1}{2}(x\wedge y)\int_0^\infty \Big[f(x+z, x+z) -
f(x,y)-f_x'(x,y)z\cr
\ar\ar\qquad\qquad\qquad\quad-f_y'(x,y)(x+z-y)\Big]\mu_{-(x - y)}(\d z)\cr
\ar\ar + \frac{1}{2}(x\wedge y)\int_0^\infty \Big[f(x+z, 2y - x + z) -
f(x,y)-f_x'(x,y)z\cr
\ar\ar\qquad\qquad\qquad\quad-f_y'(x,y)(y-x+z)\Big]\mu_{(x - y)}(\d z)\cr
\ar\ar + (x\wedge y)\int_0^\infty \Big[f(x+z, y+z) - f(x,y)- f_x'(x, y)z - f_y'(x, y)z\Big]\cr
\ar\ar\qquad\qquad\qquad\times \Big[\mu(\d z) - \frac12 \mu_{-(x - y)}(\d z) - \frac12  \mu_{(x - y)}(\d z)\Big]\cr
\ar\ar + (x - y)^+\int_0^\infty \Big[f(x + z, y) - f(x, y) - f'_x(x, y)z\Big] \mu(\d z)\cr
\ar\ar + (x - y)^-\int_0^\infty \Big[f(x, y+z) - f(x, y) - f'_y(x, y)z\Big] \mu(\d z)\cr
\ar=\ar \gamma(x)f_x'(x, y) + \gamma(y)f_y'(x, y)\cr
\ar\ar + \frac{1}{2}\sigma^2 xf_{xx}''(x, y) + \frac{1}{2}\sigma^2 yf_{yy}''(x, y)  -\sigma^2\sqrt{xy}f_{xy}''(x, y)\cr
\ar\ar + \frac{1}{2}(x\wedge y)\int_0^\infty \Big[f(x+z, x+z) - f(x+z, y+z)
\Big]\mu_{-(x - y)}(\d z)\cr
\ar\ar + \frac{1}{2}(x\wedge y)\int_0^\infty \Big[f(x+z, 2y - x + z) - f(x + z, y + z)
\Big]\mu_{(x - y)}(\d z)\cr
\ar\ar + (x\wedge y)\int_0^\infty \Big[f(x+z, y+z) - f(x,y)- f_x'(x, y)z - f_y'(x, y)z\Big]\mu(\d z)\cr
\ar\ar + (x - y)^+\int_0^\infty \Big[f(x + z, y) - f(x, y) - f'_x(x, y)z\Big] \mu(\d z)\cr
\ar\ar + (x - y)^-\int_0^\infty \Big[f(x, y+z) - f(x, y) - f'_y(x, y)z\Big] \mu(\d z),
\eeqlb
\beqlb\label{tildeL1}
\tilde{L}_ef(x, y) \ar=\ar \beta_0 x f_x'(x, y) + \beta_0 y f_y'(x, y) \cr
\ar\ar + \frac{\beta_1^2}{2}x^2f_{xx}''(x, y) + \frac{\beta_1^2}{2} y^2 f_{yy}''(x, y) -  \beta_1^2 xyf_{xy}''(x, y)
\cr
\ar\ar +\int_{\mbb{R}}\Big[f(x\e^z,y\e^z)-f(x,y)
-f_x'(x,y)(\e^z-1)x\cr
\ar\ar\qquad\quad -f_y'(x,y)(\e^z-1)y\Big] \nu(\d z)
\eeqlb
and
\beqlb\label{tildeL_c}
\tilde{L}_cf(x,y)
\ar=\ar (r(x)\wedge r(y))\int_{0}^1 \Big[f(z x, z y) -f(x,y)\Big]q(\d z)\cr
\ar\ar+(r(x)-r(y))^+\int_{0}^1 \Big[f(z x,y) -f(x,y)\Big]q(\d z)\cr
\ar\ar +(r(x)-r(y))^-\int_{0}^1 \Big[f(x,z y) -f(x,y)\Big]q(\d z).
\eeqlb
Here and in what follows, $f_x'(x, y) = \frac{\partial f(x, y)}{\partial x}, f_{xx}''(x, y) = \frac{\partial^2 f(x, y)}{\partial x^2}$ and $f_{xy}''(x, y) = \frac{\partial^2 f(x, y)}{\partial x\partial y}$.
The second equality in \eqref{tildeL0} uses the fact that $\mu_{z}(\mbb{R}_+)=\mu_{-z}(\mbb{R}_+)$ for all $z\in \mbb{R}$.
Let $\mcr{D}(\tilde{L})$ denote the linear space consisting of the functions $f$ such that the integrals in \eqref{tildeL0}, \eqref{tildeL1} and \eqref{tildeL_c} are convergent and define functions locally bounded on compact subsets of $\Delta^c.$ The operator $\tilde{L}$ determines the movement of the coupling process before
two marginal processes meet together.
Then it is easy to see that
\begin{lem} The operator $\tilde {L}$ defined by \eqref{tildeL} is a coupling operator of the generator $L$ given by \eqref{L}.\end{lem}

We call $(\tilde{L}, \mcr{D}(\tilde{L}))$ the {\it coupling generator} of the process $(X_t)_{t \ge 0},$ which is the unique strong solution to \eqref{X}.
In the following, we will construct the process $(X_t, Y_t)_{t \ge 0}$ on $\mbb{R}_+^2$ corresponding to the coupling operator $\tilde{L}$ defined \eqref{tildeL}.
Let
\beqnn
\rho(x, z) = \frac{\mu_x(\d z)}{\mu(\d z)} \in [0, 1], \qquad x \in \mbb{R}, z \in \mbb{R}_+
\eeqnn
with $\rho(0, z) = 1$ by convention.  Consider the following SDE:
\beqlb\label{XY}\left\{
\begin{aligned}
 & X_t =x + \int_0^t \gamma(X_s)\d s +\sigma\int_0^t\sqrt{X_s}\d W_s + \int_0^t\int_0^\infty\int_0^{X_{s-}} z\tilde{M}(\d s, \d z, \d u)\\
& \qquad +  \beta_0 \int_0^t X_s\d s + \beta_1\int_0^{t} X_s \d B_s   +  \int_0^{t } \int_{\mbb{R}} X_{s-}(\e^z-1) \tilde{N}(\d s, \d z)\\
& \qquad + \int_0^{t } \int_{0}^1\int_0^{r(X_{s-})} X_{s-}(z-1) Q(\d s, \d z, \d u), \\
&  Y_t = y + \int_0^{t }\gamma(Y_s)\d s + \sigma\int_0^t\sqrt{Y_s}\d W_s^* + \int_0^t\int_0^\infty\int_0^{Y_{s-}} z\tilde{M}(\d s, \d z, \d u)\cr
& \qquad + \beta_0\int_0^t Y_s \d s + \beta_1 \int_0^{t}Y_s  \d B_s^*  +  \int_0^{t } \int_{\mbb{R}} Y_{s-}(\e^z-1) \tilde{N}(\d s, \d z)\\
& \qquad + \int_0^{t } \int_{0}^1\int_0^{r(Y_{s-})} Y_{s-}(z-1) Q(\d s, \d z, \d u)  + \eta_t,
\end{aligned}
\right.
\eeqlb
where
\beqlb\label{etat}
\eta_t \ar=\ar \int_0^t U_{s-}\int_0^\infty\int_0^{\frac12(X_{s-}\wedge Y_{s-})\rho(-U_{s-}, z)} M(\d s, \d z, \d u) \cr
\ar\ar - \int_0^t U_{s-}\int_0^\infty\int_{\frac12(X_{s-}\wedge Y_{s-})\rho(-U_{s-}, z)}^{\frac12(X_{s-}\wedge Y_{s-})[\rho(-U_{s-}, z)+\rho(U_{s-}, z)]} M(\d s, \d z, \d u)
\eeqlb
with $U_t = X_t - Y_t,$ and
\beqnn
W_t^* = \left\{
\begin{aligned}
	&	-W_t, \qquad\qquad\quad t \le T,\\
	&	-2W_T + W_t, \qquad t > T,
\end{aligned}
\right.
\eeqnn
as well as
\beqnn
B_t^* = \left\{
\begin{aligned}
	&	-B_t, \qquad\qquad\quad t \le T,\\
	&	-2B_T + B_t, \qquad t > T
\end{aligned}
\right.
\eeqnn
with $T  =\inf\{t \ge 0: X_t = Y_t \}$.

\begin{prop}
There is a pathwise unique strong solution $(X_t, Y_t)_{t \ge 0}$ to the system \eqref{XY}. Moreover, it holds that $X_{T+t} = Y_{T+t} $ for every $t \ge 0$ if $T < \infty$.
\end{prop}
\proof
We first notice that there is a pathwise unique solution to the following equation:
\beqnn
Z_0(t) \ar=\ar y + \int_0^t \gamma(Z_0(s))\d s - \sigma\int_0^t\sqrt{Z_0(s)}\d W_s+ \int_0^t\int_0^\infty\int_0^{Z_0(s-)} z\tilde{M}(\d s, \d z, \d u)\cr
\ar\ar +  \beta_0 \int_0^tZ_0(s) \d s - \beta_1\int_0^{t} Z_0(s) \d B_s +  \int_0^{t } \int_{\mbb{R}} Z_0(s-)(\e^z-1) \tilde{N}(\d s, \d z)\cr
\ar\ar + \int_0^{t } \int_{0}^1\int_0^{r(Z_0(s-))} Z_0(s-)(z-1) Q(\d s, \d z, \d u).
\eeqnn
Define $T_1=\inf\{t \ge 0: X_t = Z_0(t) \}$ and $U^0(t)=X_{t} - Z_0(t)$. Let
$$\sigma^{(1)}_1  = \inf\Big\{t\ge0:  \int_0^{t\wedge T_1}\int_0^\infty\int_0^{\frac12 (X_{s-}\wedge Z_0(s-))\rho(-U^0(s-), z)}   M(\d s, \d z, \d u) =1\Big\},$$
$$\sigma^{(2)}_1  = \inf\Big\{t\ge0:  \int_0^{t\wedge T_1}\int_0^\infty\int_{\frac12(X_{s-}\wedge Z_0(s-))\rho(-U^0(s-), z)}^{\frac12(X_{s-}\wedge Z_0(s-))[\rho(-U^0(s-),  z)+\rho(U^0(s-), z)]}  M(\d s, \d z, \d u) = 1\Big\}.$$
For $t\ge 0$, let $Y_0(t)=Z_0(t\wedge T_1\wedge  \sigma^{(1)}_1\wedge  \sigma^{(2)}_1) + \eta_0(t\wedge T_1\wedge  \sigma^{(1)}_1\wedge  \sigma^{(2)}_1)$, where
\beqnn
\eta_0(t) \ar=\ar \int_0^{t\wedge T_1} U^0(s-)\int_0^\infty\int_0^{\frac12(X_{s-}\wedge Z_0(s-))\rho(-U^0(s-), z)} M(\d s, \d z, \d u) \cr
\ar\ar - \int_0^{t\wedge T_1} U^0(s-)\int_0^\infty\int_{\frac12(X_{s-}\wedge Z_0(s-))\rho(-U^0(s-), z)}^{\frac12(X_{s-}\wedge Z_0(s-))[\rho(-U^0(s-), z)+\rho(U^0(s-), z)]} M(\d s, \d z, \d u).
\eeqnn
Now we consider separately the cases $  \sigma^{(1)}_1=\sigma^{(2)}_1=\infty$, $\sigma^{(1)}_1<\sigma^{(2)}_1$ and $\sigma^{(2)}_1<\sigma^{(1)}_1$.

(i)~$\sigma^{(1)}_1=\sigma^{(2)}_1=\infty$.

In this case, we have $\eta_0 (t) = 0$ for all $t \ge 0$. Then, the process
$(X_t,Y_t)_{t \ge 0}$ is defined
by $Y_t = Z_0(t\wedge T)  + X_t - X_{t\wedge T}$, $T = T_1$ and
$\eta_t = 0$.

(ii) $\sigma^{(1)}_1<\sigma^{(2)}_1\le \infty$.

In this case, we have $\sigma^{(1)}_1 \le T_1$. Moreover, we can define $Y_t = Z_0 (t)$ and $\eta_t = 0$ for $0 \le t <\sigma^{(1)}_1$. If $U^0(\sigma^{(1)}_1-)\ge0$, then
\beqnn
Y_0(\sigma^{(1)}_1)\ar: = \ar Z_0(\sigma^{(1)}_1) + \eta_0(\sigma^{(1)}_1) = Z_0(\sigma^{(1)}_1) + \Delta \eta_0(\sigma^{(1)}_1)\\
\ar = \ar   Z_0(\sigma^{(1)}_1) + U^0(\sigma^{(1)}_1-)  \ge 0.
\eeqnn
If $U^0(\sigma^{(1)}_1-)<0$, then similar to the proof of Lemma 3.3 in \cite{LLWZ23} by using $\rho(x, z) = 0$ for $0 < z \le 0 \vee x$,
\beqnn
\Delta Z_0(\sigma^{(1)}_1) \ar\ar=\int_{\{\sigma^{(1)}_1\}}\int_0^\infty\int_{0}^{\frac12\big(X_{\sigma^{(1)}_1-}\wedge Z_0(\sigma^{(1)}_1-)\big) \rho (- U^0(\sigma^{(1)}_1-), z ) } z M(\d s, \d z, \d u)\\
\ar\ar >- U^0(\sigma^{(1)}_1-)
\eeqnn
and
\beqnn
Y_0(\sigma^{(1)}_1)
\ar := \ar   Z_0(\sigma^{(1)}_1-) +\Delta Z_0(\sigma^{(1)}_1) +  U^0(\sigma^{(1)}_1-) \ge 0.
\eeqnn

(iii) $\sigma^{(2)}_1<\sigma^{(1)}_1\le \infty$.

In this case, we have $\sigma^{(2)}_1 \le T_1$. Moreover, we can define $Y_t = Z_0 (t)$ and $\eta_t = 0$ for $0 \le t <\sigma^{(2)}_1$. If $U^0(\sigma^{(2)}_1-)\le 0$, then
\beqnn
Y_0(\sigma^{(2)}_1)\ar := \ar Z_0(\sigma^{(2)}_1) + \eta_0(\sigma^{(2)}_1) = Z_0(\sigma^{(2)}_1) + \Delta \eta_0(\sigma^{(2)}_1)\\
\ar = \ar   Z_0(\sigma^{(2)}_1) - U^0(\sigma^{(2)}_1-)  \ge0.
\eeqnn
If $U^0(\sigma^{(2)}_1-)>0$, then by using $\rho(x, z) = 0$ for $0 < z \le 0 \vee x$ again,
\beqnn
 \Delta Z_0(\sigma^{(2)}_1) \ar=\ar\int_{\{\sigma^{(2)}_1\}}\int_0^\infty\int_{\frac12\big(X_{\sigma^{(2)}_1-}\wedge Z_0(\sigma^{(2)}_1-)\big)\rho (- U^0(\sigma^{(2)}_1-), z ) }^{\frac12\big(X_{\sigma^{(2)}_1-}\wedge Z_0(\sigma^{(2)}_1-)\big) [\rho (- U^0(\sigma^{(2)}_1-), z ) +\rho ( U^0(\sigma^{(2)}_1-), z ) ]} z M(\d s, \d z, \d u)\\
\ar>\ar U^0(\sigma^{(2)}_1-)
\eeqnn
and
\beqnn
Y_0(\sigma^{(2)}_1)
\ar: = \ar   Z_0(\sigma^{(2)}_1-) +\Delta Z_0(\sigma^{(2)}_1) +  U^0(\sigma^{(2)}_1-)  \ge 0.
\eeqnn

Below we consider the construction for $t\ge \sigma^{(i)}_1,i=1,2$. Let $X_1 (t) = X_{\sigma^{(i)}_1+t} $ for $t \ge 0$. From \eqref{X} it follows that
\beqnn
X_1 (t) \ar=\ar X_{\sigma^{(i)}_1}+ \int_{0}^{t} \gamma(X_1(s))\d s+ \sigma\int_{0}^{t}\sqrt{X_1(s)}\d W_{s+\sigma^{(i)}_1}\cr
\ar\ar  + \int_{0}^{t}\int_0^\infty\int_0^{X_1(s-)} z\tilde{M}(\sigma^{(i)}_1+\d s, \d z, \d u)\cr
\ar\ar + \beta_0 \int_{0}^{t} X_1(s) \d s + \beta_1 \int_{0}^{t} X_1(s)\d B_{s+\sigma^{(i)}_1}\cr
\ar\ar+ \int_{0}^{t}X_1(s-)\int_{\mbb{R}}(\e^z - 1) \tilde{N}(\sigma^{(i)}_1+\d s, \d z)\cr
\ar\ar+ \int_{0}^{t}X_1(s-)\int_{0}^1\int_0^{r(X_1(s-))} (z - 1) Q(\sigma^{(i)}_1+\d s, \d z,\d u).\cr
\eeqnn
We can also construct $(Z_1 (t))_{t \ge 0}$ by the pathwise unique solution to
\beqnn
Z_1(t)\ar=\ar Y_0(\sigma^{(i)}_1) + \int_0^t \gamma(Z_1(s))\d s + \int_0^t\sqrt{Z_1(s)}\d W_{s+\sigma^{(i)}_1}^*\cr
\ar\ar + \int_0^t\int_0^\infty\int_0^{Z_1(s-)} z\tilde{M}(\sigma^{(i)}_1+\d s, \d z, \d u)\cr
\ar\ar + \beta_0 \int_0^t Z_1(s) \d s + \beta_1 \int_0^{t}Z_1(s)  \d B_{s+\sigma^{(i)}_1}^*\cr \ar\ar +  \int_0^{t } \int_{\mbb{R}} Z_1(s-)(e^z-1) \tilde{N}(\sigma^{(i)}_1+\d s, \d z)\cr
\ar\ar + \int_0^{t } \int_{0}^1\int_0^{r(Z_1(s-))} Z_1(s-)(z-1) Q(\sigma^{(i)}_1+\d s, \d z, \d u).\cr
\eeqnn
Then we repeat the procedure of $(X_t, Z_0(t))_{t \ge 0}$  for
the process $(X_1(t), Z_1(t))_{t \ge 0}$.

Similar to \cite{LW20},  only finitely many modifications have to be made in the interval $(0,t \wedge\tau_m)$, where
$\tau_m = \inf\{t \ge 0 : Y_t > m$  or  $|X_t - Y_t | < 1/m\}.$
Finally, by letting $m \rightarrow\infty$, we obtain the unique strong solution to the SDE \eqref{XY}
globally. The second assertion is a direct consequence of the construction above.
\qed

Finally, we have the following statement.

\begin{cor}The unique strong solution $(X_t, Y_t)_{t \ge 0}$ to the system \eqref{XY} is a Markov coupling of the process $(X_t)_{t\ge0}$ determined by \eqref{X}.\end{cor}
\proof By the It\^{o}'s formula, one can get that the infinitesimal generator of the unique strong solution $(X_t, Y_t)_{t \ge 0}$ to the system \eqref{XY} is just the coupling operator $\tilde {L}$ defined by \eqref{tildeL}. Then, due to the uniqueness of the strong solution to the  system \eqref{XY}, we can obtain the desired assertion. \qed

\section{General result for the exponential ergodicity}

In this section, a general result for the exponential ergodicity of the process $(X_t)_{t \ge 0}$ will be given.
We raise the following four conditions before giving the main result. We call a $C^{2}$-function $V\in {\mcr D}(L)$ a {\it Lyapunov function} for the process $(X_t)_{t\ge0}$, if $V\ge1$, and there are constants $\lambda_1>0$ and $\lambda_2>0$ such that
\beqlb\label{lyapunov}
LV(x)\leq\lambda_2-\lambda_1 V(x),\quad\quad x\in\mbb{R}_+,
\eeqlb
where $L$ is the infinitesimal generator given by (\ref{L}).

\bcond\label{Ncond}
The rate function $r(x)$ for the catastrophes is globally Lipschitz, and the constant $\alpha=\gamma(0)>0$ in the drift term $\gamma(x)$.
\econd

\bcond\label{Lcond}
{\bf(Lyapunov condition)} There exists a Lyapunov function $V(x)$ satisfying $V(x)\rightarrow\infty$ as $x\rightarrow\infty$.
\econd

\bcond\label{Fcond}
{\bf(Non-triviality of branching mechanism)} We have either $\sigma^2 > 0$ or that $\int_0^1 z\mu(\d z)=\infty$ and there exist constants $c_0 > 0$ and $\delta > 0$ so that
for all $|x|\le c_0$,
\beqnn
\mu_x(\mbb{R}_+)\ge \delta.
\eeqnn
\econd

\bcond\label{Ccond}
\beqnn
\limsup\limits_{x\rightarrow \infty}\frac{H(x)}{V(x)}=0,
\eeqnn
where
\beqnn
H(x):= \int_{0}^{1/x} (1 -  zx  )^3\big( \nu(\d \ln z)+  r(x)q(\d z)\big).
\eeqnn
\econd

The main theorem in this section is as follows.

\bth\label{main result}
Suppose that Condition \ref{Ncond}--Condition \ref{Ccond} are satisfied. Then the solution to \eqref{X} is exponentially ergodic in the $V$-weighted total variation distance.
\eth

 Theorem \ref{main result} is   more general than Theorem \ref{example}. The existence of a suitable Lyapunov function has become a standard condition for the ergodicity of Markov processes; see, e.g., \cite{LLWZ23,LMW21,Mas07,MT92,MT93a,MT93b}. In particular, from Theorem \ref{example} and its proof below one can see that Condition \ref{Lcond} roughly indicates that all the competition mechanism, catastrophes and environments could help to guarantee the exponential ergodicity of the process.

To prove Theorem \ref{main result}, the main task is to find a distance-like function $F\in \mcr{D}(\tilde{L})$ such that
\beqnn
\tilde{L}F(x, y) \le -\lambda F(x, y), \qquad (x, y) \in \Delta^c,
\eeqnn where $\tilde{L}$ is the coupling operator constructed in Subsection \ref{subsection2.2}.
For this, we will make full use of the Lyapunov function for an unbounded area, and utilize inner structure of the process reflected by the coupling generator for the bounded area.

\subsection{Estimate of the coupling generator}\label{section estimate}
Recall that $\mcr{D}(\tilde{L})$ is the linear space consisting of the functions $f$ such that the integrals in \eqref{tildeL0}, \eqref{tildeL1} and \eqref{tildeL_c} are convergent and define functions locally bounded on compact subsets of $\Delta^c$, and that $C_b^2(\mbb{R}_+)$ denotes the space of bounded and continuous functions on $\mbb{R}_+$ with bounded and continuous derivatives up to the second order. For any $l_0>0$, define
\beqlb\label{fxy}
f(x, y) =\phi(x\vee y) \psi(|x - y|\wedge l_0)1_{\{x \neq y\}},
\eeqlb
where $\psi \in C_b^2(\mbb{R}_+)$ is a nonnegative and concave nondecreasing function with $\psi(0)=1$ and  $\phi\in C_b^2(\mbb{R}_+)$ is a nonnegative nonincreasing function.
In particular, $f(z, z) = 0$ for any $z \ge 0,$  and, by \eqref{tildeL}-\eqref{tildeL_c}, one sees that $f \in \mcr{D}(\tilde{L}).$

\subsubsection{Preliminary estimation of the coupling generator}\label{section preliminary estimate}
In the following, we give the preliminary estimation of $\tilde{L}f(x,y)$ for $x>y$ according to different structures of the CBIRE-processes with competition and catastrophes. The case for $y>x$ can be discussed similarly.

(i) ${Branching}$

Recall that $\gamma(x) = \alpha - bx- g(x)$, where $g$ is nondecreasing. Moreover, $\psi$ is nondecreasing and concave on $\mbb{R}_+$, which implies that $\psi(2r) - 2\psi(r) \le -\psi(0)=-1$ for all $r\ge0$. Then, by \eqref{tildeL0}, for $0<x-y\le l_0$,
\beqnn
\tilde{L}_bf(x, y) \ar=\ar { \gamma(x)\phi'(x)\psi(x - y) + (\gamma(x) - \gamma(y))\phi(x)\psi'(x - y)}\cr
\ar\ar + {  \frac12 \sigma^2 x\phi''(x)\psi(x - y)} + \frac12 \sigma^2(x + y) \phi(x)\psi''(x - y) + \sigma^2x\phi'(x)\psi'(x - y)\cr
\ar\ar + \sigma^2\sqrt{xy}\phi'(x)\psi'(x - y) + \sigma^2\sqrt{xy}\phi(x)\psi''(x - y) \cr
\ar\ar - \frac{1}{2}y\psi(x - y)\int_0^\infty \phi(x + z)  \mu_{-(x - y)}(\d z)\cr
\ar\ar + {  \frac12 y [\psi( (2(x - y))\wedge l_0) -\psi(x - y)]\int_0^\infty \phi(x+ z)\mu_{(x - y)}(\d z)}\cr
\ar\ar +  { x\psi(x - y)\int_0^\infty [\phi(x + z) - \phi(x) - z\phi'(x)] \mu(\d z)}\cr
\ar\ar + (x - y)\int_0^\infty \left[\phi(x + z)\left(\psi((x - y +z)\wedge l_0) - \psi(x - y)\right) - z\phi(x)\psi'(x - y)\right] \mu(\d z)\cr
\ar\le\ar {  \gamma(x)\phi'(x)\psi(x - y) -b (x - y)\phi(x)\psi'(x - y)}\cr
\ar\ar + { \frac12 \sigma^2 x\phi''(x)\psi(x - y)} + \frac12 \sigma^2(x + y) \phi(x)\psi''(x - y)\cr
\ar\ar {  - \frac{1}{2}y\psi(x - y)\int_0^\infty [\phi(x + z)-\phi(x)]  \mu_{-(x - y)}(\d z)}- \frac12y\phi(x)\mu_{(x - y)}(\mbb{R}_+) \cr
\ar\ar +  {  x\psi(x - y)\int_0^\infty [\phi(x + z) - \phi(x) - z\phi'(x)] \mu(\d z)}\cr
\ar\ar + (x - y)\phi(x)\int_0^\infty \left[\psi((x - y +z)\wedge l_0) - \psi(x - y) - z\psi'(x - y)\right] \mu(\d z);
\eeqnn
and for $x-y>l_0$,
\beqnn
\tilde{L}_bf(x, y) \ar=\ar { \gamma(x)\phi'(x)\psi(l_0) + \frac{1}{2}\sigma^2 x\phi''(x)\psi(l_0) }\cr
\ar\ar - \frac{1}{2}y\psi(l_0)\int_0^\infty \phi(x + z)  \mu_{-(x - y)}(\d z)\cr
\ar\ar +{  x\psi(l_0)\int_0^\infty [\phi(x + z) - \phi(x) - z\phi'(x)] \mu(\d z)}.
\eeqnn
Here we note that, in the case $y>x$,  the fifth term in the equality above  for $0<x-y\le l_0$ turns to be  \begin{align*}&\frac12 x \int_0^\infty [\phi(2y-x+ z)\psi(l_0 \wedge 2(y - x)) - \phi(y + z)\psi(y - x)]\mu_{(x - y)}(\d z)\\
&\le \frac12 x [\psi( (2(y - x))\wedge l_0) -\psi(x - y)]\int_0^\infty \phi(y+ z)\mu_{(x - y)}(\d z) \end{align*} since $2y-x+ z>y+z$ and $\phi(2y-x+ z)\le \phi(y + z)$. Hence in this case the equality here should become an inequality.

(ii) ${Random~environment}$

By \eqref{tildeL1}, for $0<x-y\le l_0$,
\beqnn
\tilde{L}_ef(x, y)
\ar=\ar \beta_0 x \phi'(x)\psi(x-y) +{ \beta_0(x - y)\phi(x)\psi'(x- y)}\cr
\ar\ar + {\frac{\beta_1^2}{2}x^2\phi''(x)\psi(x - y)} + \frac{\beta_1^2}{2}(x^2+y^2) \phi(x)\psi''(x - y) + \beta_1^2x^2\phi'(x)\psi'(x - y)\cr
\ar\ar + \beta_1^2xy\phi'(x)\psi'(x - y) + \beta_1^2xy\phi(x)\psi''(x - y) \cr
\ar\ar + { \int_{\mbb{R}}\Big[\phi(xe^z) \psi(l_0 \wedge (\e^z(x - y))) - \phi(x)\psi(x - y)}\cr
\ar\ar\qquad\quad {- (\e^z - 1)\left[ x \phi'(x)\psi(x - y) +  (x - y)\phi(x)\psi'(x - y)\right]\Big]\nu(\d z)}\cr
\ar\le \ar \beta_0 x \phi'(x)\psi(x-y) +{  \beta_0(x - y)\phi(x)\psi'(x- y)}\cr
\ar\ar + {  \frac{\beta_1^2}{2}x^2\phi''(x)\psi(x - y)} + \frac{\beta_1^2}{2}(x+y)^2 \phi(x)\psi''(x - y) \cr
\ar\ar + \phi(x)\int_{\mbb{R}}\left[\psi((\e^z(x - y))\wedge l_0) - \psi(x-y) - (\e^z - 1)(x-y)\psi'(x - y)\right]\nu(\d z)\cr
\ar\ar + \int_{\mbb{R}} \left[\phi(xe^z) - \phi(x)\right]\left[\psi((\e^z(x - y))\wedge l_0) - \psi(x - y)\right]\nu(\d z)\cr
\ar\ar + { \psi(x - y)\int_{\mbb{R}}\left[\phi(x\e^z) - \phi(x) - (\e^z - 1)x \phi'(x)\right]\nu(\d z)};
\eeqnn
and for $x-y>l_0$,
\beqnn
\tilde{L}_ef(x, y) \ar=\ar {  \beta_0 x\phi'(x)\psi(l_0) + \frac{\beta_1^2}{2}x^2\phi''(x)\psi(l_0) } \cr
\ar\ar + \phi(x)\int_{\mbb{R}}\left[\psi((\e^z(x - y))\wedge l_0) - \psi(l_0) \right]\nu(\d z)\cr
\ar\ar + \int_{\mbb{R}} \left[\phi(x\e^z) - \phi(x)\right]\left[\psi((\e^z(x - y))\wedge l_0) - \psi(l_0)\right]\nu(\d z)\cr
\ar\ar + {  \psi(l_0)\int_{\mbb{R}}\left[\phi(x\e^z) - \phi(x) - (\e^z - 1)x \phi'(x)\right]\nu(\d z)}.
\eeqnn

(iii) ${Catastrophes}$

By \eqref{tildeL_c}, for $0<x-y\le l_0$,
\beqnn
\tilde{L}_cf(x, y) \ar=\ar \left(r(x)\wedge r(y)\right)\int_0^1 \left[\phi(zx)\psi(z(x - y)) - \phi(x)\psi(x-y)\right]q(\d z)\cr
\ar\ar + (r(x) - r(y))^+\int_0^1\left[\phi(zx\vee y)\psi(|zx-y|\wedge l_0) -\phi(x)\psi(x-y)\right]q(\d z)\cr
\ar\ar + (r(x) - r(y))^-\int_0^1\phi(x)\left[\psi((x-zy)\wedge l_0) - \psi(x - y)\right]q(\d z)\cr
\ar=\ar  \left(r(x)\wedge r(y)\right)\int_0^1 \psi(z(x - y))\left[\phi(zx) - \phi(x)\right]q(\d z)\cr
\ar\ar +  \left(r(x)\wedge r(y)\right)\int_0^1 \phi(x)\left[\psi(z(x - y)) - \psi(x - y)\right]q(\d z)\cr
\ar\ar+ (r(x) - r(y))^-\int_0^1\phi(x)\left[\psi((x-zy)\wedge l_0) - \psi(x - y)\right]q(\d z)\cr
\ar\ar + (r(x) - r(y))^+\int_0^1\phi((zx)\vee y)[\psi(|zx-y|\wedge l_0) -\psi(x-y)]q(\d z)\cr
\ar\ar+ (r(x) - r(y))^+\psi(x-y)\int_0^1\left[\phi((zx)\vee y) -\phi(x)\right]q(\d z)\cr
\ar\le\ar {  \left(r(x)\wedge r(y)\right) \psi(x - y)\int_0^1\left[\phi(zx) - \phi(x)\right]q(\d z)}\cr
\ar\ar- \left(r(x)\wedge r(y)\right) \phi(x)\psi'(x - y)(x - y)\int_0^1 (1 - z)q(\d z)\cr
\ar\ar + {  (r(x) - r(y))^-\psi'(x - y)(l_0 - (x - y))\int_0^1\phi(x)q(\d z)}\cr
\ar\ar +{   (r(x) - r(y))^+\psi'(x - y)(l_0 - (x - y))\int_0^1 \phi(zx)q(\d z)}\cr
\ar\ar + {  (r(x) - r(y))^+ \psi(x - y)\int_0^1\left[\phi(zx) - \phi(x)\right]q(\d z)}\cr
\ar  \le \ar{ r(x)\psi(x - y)\int_0^1[\phi(zx) - \phi(x)]q(\d z)}\cr
\ar\ar + { |r(x) - r(y)|\psi'(x - y)(l_0 - (x - y))\int_0^1\phi(zx)q(\d z)}\cr
\ar\ar - (r(x)\wedge r(y))\phi(x)\psi'(x - y)(x - y)\int_0^1(1-z )q(\d z);
\eeqnn
and for $x-y>l_0$,
\beqnn
\tilde{L}_cf(x,y)
\ar=\ar (r(x)\wedge r(y))\int_0^1 \left[\phi(zx)\psi((z(x - y))\wedge l_0) - \phi(x)\psi(l_0)\right]q(\d z)\cr
\ar\ar + (r(x) - r(y))^+\int_0^1\phi((zx)\vee y)[\psi(|zx-y|\wedge l_0) -\psi(l_0)]q(\d z)\cr
\ar\ar+ (r(x) - r(y))^+\psi(l_0)\int_0^1\left[\phi((zx)\vee y) -\phi(x)\right]q(\d z)\cr
\ar\le\ar (r(x)\wedge r(y))\int_0^1 \phi(zx)\left[\psi((z(x - y))\wedge l_0) - \psi(l_0)\right]q(\d z)\cr
\ar\ar+{ r(x)\psi(l_0)\int_0^1 \left[\phi(zx)  - \phi(x)\right]q(\d z)}.
\eeqnn
In conclusion, we have for $0<x-y\le l_0$,
\beqlb \label{L_dlel}
\tilde{L}f(x,y)\ar\le\ar  \psi(x - y)\Big[{ (\gamma(x)+\beta_0 x) \phi'(x)}+\frac12 x(\sigma^2 + x\beta_1^2) \phi''(x)\cr
\ar\ar\qquad\qquad\quad {  - \frac{1}{2}y\int_0^\infty [\phi(x + z)-\phi(x)]  \mu_{-(x - y)}(\d z)}\cr
\ar\ar\qquad\qquad\quad +  {  x\int_0^\infty [\phi(x + z) - \phi(x) - z\phi'(x)] \mu(\d z)}\cr
\ar\ar\qquad\qquad\quad +{  \int_{\mbb{R}}\left[\phi(x\e^z) - \phi(x) - (\e^z - 1)x \phi'(x)\right]\nu(\d z)}\cr
\ar\ar\qquad\qquad\quad +{   r(x)\int_0^1[\phi(zx) - \phi(x)]q(\d z)}\Big]\cr
 \ar\ar + \phi(x)\Big\{(x - y) \Big[{ (\beta_0 -b)\psi'(x - y)}+ \frac12\sigma^2\psi''(x - y)\cr
\ar\ar\qquad\qquad\qquad\qquad+\int_0^\infty \left[\psi(x - y +z) - \psi(x - y) - z\psi'(x - y)\right] \mu(\d z)\cr
\ar\ar\qquad\qquad\qquad\qquad- (r(x)\wedge r(y))\psi'(x - y)\int_0^1(1-z )q(\d z)\Big]\cr
\ar\ar\qquad\qquad + \frac12\sigma^2 y\psi''(x - y)- \frac12y\mu_{(x - y)}(\mbb{R}_+)+\frac12\beta_1^2(x - y)^2\psi''(x - y) \cr
\ar\ar\qquad\qquad + \int_{\mbb{R}}\left[\psi(\e^z(x - y)) - \psi(x-y) - (\e^z - 1)(x-y)\psi'(x - y)\right]\nu(\d z)\Big\}\cr
\ar\ar + { |r(x) - r(y)|\psi'(x - y)l_0 \int_0^1\phi(zx)q(\d z)};
\eeqlb
and for $x-y>l_0$,
\beqlb\label{L_dgl}
\tilde{L}f(x,y)\ar\le\ar \psi(l_0)\Big[{ (\gamma(x)+\beta_0x)\phi'(x) + \frac{1}{2} x(\sigma^2+x\beta_1^2)\phi''(x)}\cr
\ar\ar\qquad\quad- \frac{1}{2}y\int_0^\infty \phi(x + z)  \mu_{-(x - y)}(\d z)\cr
\ar\ar\qquad\quad +{  x\int_0^\infty [\phi(x + z) - \phi(x) - z\phi'(x)] \mu(\d z)}\cr
\ar\ar\qquad\quad + {  \int_{\mbb{R}}\left[\phi(x\e^z) - \phi(x) - (\e^z - 1)x \phi'(x)\right]\nu(\d z)}\cr
\ar\ar\qquad\quad + {  r(x)\int_0^1 \left[\phi(zx)  - \phi(x)\right]q(\d z)}\Big] \cr
\ar\ar + \phi(x)\int_{\mbb{R}}\left[\psi((\e^z(x - y))\wedge l_0) - \psi(l_0) \right]\nu(\d z)\cr
\ar\ar + \int_{\mbb{R}} \left[\phi(x\e^z) - \phi(x)\right]\left[\psi((\e^z(x - y))\wedge l_0) - \psi(l_0)\right]\nu(\d z)\cr
\ar\ar +(r(x)\wedge r(y))\int_0^1 \phi(zx)\left[\psi((z(x - y))\wedge l_0) - \psi(l_0)\right]q(\d z).
\eeqlb
\subsubsection{Detailed estimation of the coupling generator}\label{section detailed estimate}
In the previous section, we give a preliminary estimation for the coupling generator $\tilde{L}$ acting on the function $f$ defined by \eqref{fxy}. To move further, we should take the especial form of the function $f$ by taking explicit $\psi$ and $\phi$ in \eqref{fxy}. In this part, we still consider $x>y$ only and the case that $y>x$ can be obtained in the similar manner.

Let $K>0$ and $ x_0\in (0,1\wedge c_0]$, and set $K_0=\min\{K,{6\alpha}/{x_0}\}$, where
$c_0>0$ is given in Condition \ref{Fcond} and
 $\alpha>0$ is the constant in the drift term $\gamma(x)$.
Recall that $V(x)$ and $H(x)$ are functions given in Conditions \ref{Lcond} and \ref{Ccond} respectively. By Condition \ref{Ccond},
\beqnn
  \liminf \limits_{x\rightarrow \infty}\Big(1-\frac{9\lambda_2}{K_0\lambda_1}\frac{H(x)}{V(x)}\Big)=1,
\eeqnn where $\lambda_1$ and $\lambda_2$ are given in Condition \ref{Lcond}. Moreover, by
$V\in C^2(\mbb{R}_+)$ and $V(x)\to \infty$ as $x\to \infty$,  there exists $M:=M(\lambda_1,\lambda_2,K_0)\ge1$ such that for $x\ge M$,
\beqnn
V(x)\ge 12~\mathrm{and}~1-\frac{9\lambda_2}{K_0\lambda_1}\frac{H(x)}{V(x)}\ge \frac{1}{4}.
\eeqnn
Let
\beqlb
S_0 \ar=\ar \{(x,y): \lambda_1(V(x)+V(y))\leq { 6\lambda_2}\},\\
l_0 \ar=\ar \sup_{(x,y)\in S_0}(|x-y|)+M.
\eeqlb

Recall that $\sigma^2 > 0$ or $\int_0^1 z\mu(\d z)=\infty$ in Condition \ref{Fcond}, one sees that
\begin{equation}\label{e:need}
\frac{\Phi(\lambda)}{\lambda}=b+\sigma^2\lambda+\int_0^\infty\frac{\e^{-\lambda z}-1+\lambda z}{\lambda}\mu(\d z)\rightarrow b+\sigma^2\cdot\infty+\int_0^\infty z\mu(\d z)=\infty
\end{equation}as $\lambda\to \infty,$ where $\Phi$ is the branching mechanism given by \eqref{branching mechanism}.
Then, there exists a constant ${ \lambda_3} > 0$ such that
\beqlb\label{R}
\tilde{\Phi}({ \lambda_3}):= \Phi({ \lambda_3}) + \left[\inf_{x \ge 0}(r(x))\int_0^1(1 - z)q (\d z) - \beta_0\right] { \lambda_3} > 0.
\eeqlb
Below, for $\lambda_0 > {\lambda_3}$ (with $\lambda_3$ given in  (\ref{R})) and $\theta\ge 4$ which are related to $l_0\ge1$  to be specified later, define
\beqlb\label{psi}
\psi(x) = 2 - \e^{-\lambda_0 x}, \qquad x \ge 0,
\eeqlb
and
\beqlb\label{phi}
\phi(x) =
\begin{cases}
	\theta + (1 - x/x_0)^3, & 0 \le x < x_0,\\
	\theta, & x \ge x_0.
\end{cases}
\eeqlb
It is easy to see that $1 \le \psi(x) \le 2,$  $\psi'(x) = \lambda_0\e^{-\lambda_0x}$ and $\psi''(x) = -\lambda_0^2\e^{-\lambda_0x}$ for any $x \ge 0;$ and that
 $\theta \le \phi(x) \le \theta + 1 \le 2\theta$ for all $x\ge0$.

Due to \eqref{psi}, we have for $u=x-y>0$,
\beqnn
\ar\ar{(\beta_0 -b)\psi'(u)}+ \frac12\sigma^2\psi''(u)+\int_0^\infty \left[\psi(u +z) - \psi(u) - z\psi'(u)\right] \mu(\d z)\cr
\ar\ar\qquad- (r(x)\wedge r(y))\psi'(u)\int_0^1(1-z )q(\d z)\cr
\ar \ar= -\e^{-\lambda_0 u}\Big[\Phi(\lambda_0)+\Big((r(x)\wedge r(y))\int_0^1(1-z )q(\d z)-\beta_0\Big)\lambda_0\Big].
\eeqnn
Moreover, for $0< u=x-y\le l_0$,
\beqnn
\ar\ar\frac12\beta_1^2u^2\psi''(u)+\int_{\mbb{R}}\left[\psi(\e^zu) - \psi(u) - (\e^z - 1)u\psi'(u)\right]\nu(\d z)\cr
\ar \ar \le \frac{1}{2}\psi''(u)u^2\Big[\beta_1^2+\int_{-\infty}^0(\e^z-1)^2\nu(\d z) + \e^{-\lambda_0 (\e-1)u}\int_{0}^1(\e^z-1)^2\nu(\d z)\Big]\cr
\ar \ar\le \frac{1}{2}\psi''(u)u^2\Big[\beta_1^2+\int_{-\infty}^0(\e^z-1)^2\nu(\d z) + \e^{-\lambda_0 (\e-1)l_0}\int_{0}^1(\e^z-1)^2\nu(\d z)\Big]\cr
\ar\ar= \frac{1}{2}\psi''(u)u^2E(\lambda_0,l_0),
\eeqnn
where \beqlb
E(\lambda_0,l_0)\ar:=\ar \beta_1^2+\int_{-\infty}^0(\e^z - 1)^2 \nu(\d z) +\e^{-\lambda_0(\e - 1)l_0}\int_0^1(\e^z - 1)^2\nu(\d z)\ge 0
\eeqlb
represents the impact of the fluctuation of random environment.

As for $\phi$, we have $\phi''(x) \le \frac{6}{x_0^2}1_{\{x \le x_0\}}$, and for $z\ge 0$,
\beqnn
\ar\ar 0 \le \phi(x) - \phi(x + z) \le 1_{\{x \le x_0\}},\cr
\ar\ar 0 \le -\phi'(x)z \le \frac{3z}{x_0}1_{\{x \le x_0\}},\cr
\ar\ar \phi(x + z) - \phi(x) - \phi'(x)z\le \frac{3z}{x_0}1_{\{x \le x_0\}},
\eeqnn
and for $z\in\mbb{R}$,
\beqnn
\phi(x + z) - \phi(x) - \phi'(x)z \le \frac{3z^2}{x_0^2}1_{\{x \le x_0\}}.
\eeqnn
On the other hand, for $x> x_0$, by the fact that $\phi(x) = \theta$ for all $x>x_0$,
\beqnn
\ar\ar \int_{\mbb{R}}\left[\phi(x\e^z) - \phi(x) - (\e^z - 1)x \phi'(x)\right]\nu(\d z) +   r(x)\int_0^1[\phi(zx) - \phi(x)]q(\d z)\cr
\ar \ar =\int_{0}^{x_0/x}\left[\phi(zx) - \phi(x)\right]\nu(\d \ln z) +   r(x)\int_0^{x_0/x}[\phi(zx) - \phi(x)]q(\d z)\cr
\ar \ar= \int_{0}^{x_0/x}\left(1 - \frac{zx}{x_0}\right)^3\big( \nu(\d \ln z)+  r(x)q(\d z)\big)=:H(x,x_0).
\eeqnn
Here, $H(x,x_0)$ represents the impact of the negative jump.

The following lemma gives the estimation of  $\tilde{L}f(x, y)$ for $x>x_0$.
\blem\label{x>x0}
	For $x> x_0$ and $0 < x - y \le l_0$,
\beqlb\label{case1}
\tilde{L}f(x, y)
\ar\le\ar  2{H(x,x_0) }- \frac12 \theta y[\sigma^2 \lambda_0^2\e^{-\lambda_0 (x - y)} + \mu_{(x - y)}(\mbb{R}_+)] \cr
\ar\ar -\theta(x - y)\lambda_0\e^{-\lambda_0 (x - y)} \bigg[{\Phi(\lambda_0)/\lambda_0+\Big((r(x)\wedge r(y))\int_0^1(1-z )q(\d z)-\beta_0\Big)}\cr
\ar\ar\qquad\qquad\qquad\qquad\qquad+\frac{1}{2}\lambda_0(x - y)E(\lambda_0,l_0)
-{ 2  l_0  \frac{|r(x) - r(y)|}{x-y}  }\bigg];
\eeqlb
and for $x> x_0$ and $ x - y >l_0$, $$\tilde{L}f(x,y)  \le  2{ H(x,x_0)}.$$
\elem

\proof
By \eqref{L_dlel} and the properties of $\psi$ and $\phi$ given as above, we have for $x> x_0$ and $0 < x - y \le l_0$,
\beqnn
\tilde{L}f(x,y)\ar\le\ar  \psi(x - y)H(x,x_0)\cr
\ar\ar + \theta\Big\{ -(x - y)\e^{-\lambda_0 (x - y)}\Big[\Phi(\lambda_0)+\Big((r(x)\wedge r(y))\int_0^1(1-z )q(\d z)-\beta_0\Big)\lambda_0\Big]\cr
\ar\ar\qquad\qquad + \frac12\sigma^2 y\psi''(x - y)- \frac12y\mu_{(x - y)}(\mbb{R}_+)   + \frac{1}{2}(x - y)^2\psi''(x - y)E(\lambda_0,l_0)\Big\}\cr
\ar\ar
+ { |r(x) - r(y)|\psi'(x - y)l_0 \int_0^1\phi(zx)q(\d z)}\cr
\ar\le\ar  2{ H(x,x_0) }
+2\theta l_0 { |r(x) - r(y)|\lambda_0\e^{-\lambda_0 (x - y)} q((0,1])}\cr
\ar\ar + \theta\Big\{ -(x - y)\e^{-\lambda_0 (x - y)}\Big[\Phi(\lambda_0)+\Big((r(x)\wedge r(y))\int_0^1(1-z )q(\d z)-\beta_0\Big)\lambda_0\Big]\cr
\ar\ar\qquad\qquad - \frac12 y[\sigma^2 \lambda_0^2\e^{-\lambda_0 (x - y)} + \mu_{(x - y)}(\mbb{R}_+)]  -\frac{1}{2}(x - y)^2\lambda_0^2\e^{-\lambda_0 (x - y)}E(\lambda_0,l_0)\Big\}.
\eeqnn
On the other hand, by \eqref{L_dgl}, it is easy to see that for $x> x_0$ and $  x - y > l_0$,
\beqnn
\tilde{L}f(x,y) \le  \psi(l_0){  H(x,x_0)} \le  2{ H(x,x_0)}.
\eeqnn
The result follows.
\qed

Note that $ x_0\le 1$ and $l_0\ge 1$. Then $x\le x_0$ implies that $0 <  x - y  \le l_0$. Next we give the estimation of  $\tilde{L}f(x, y)$ for $x\le x_0$.
\blem\label{xlex0} For $x  \le x_0$,
\beqlb\label{case2}
\tilde{L}f(x, y)
\ar\le\ar { (\gamma(x)+\beta_0 x) \phi'(x) \psi(x - y)}-\frac12y[(\theta- 2)\mu_{(x - y)}(\mbb{R}_+) + \sigma^2 \theta \lambda_0^2\e^{-\lambda_0(x - y)}] \cr
\ar\ar +  {  \frac{6 x}{x_0^2}\Big[\Big(\sigma^2 +  \int_0^1 z^2 \mu(\d z)\Big) + x \Big( \beta_1^2  + \int_{-\infty}^1(\e^z - 1)^2\nu(\d z)\Big)\Big]}\cr
\ar\ar + {  \frac{6x}{x_0}\Big[\int_1^\infty z  \mu(\d z)+\int_1^\infty(\e^z - 1)\nu(\d z) + r(x)\int_0^1(1 - z)q(\d z)\Big]}\cr
\ar\ar  - \phi(x)\lambda_0\e^{-\lambda_0 (x - y)}(x - y)\bigg[  { \Phi(\lambda_0)/\lambda_0+\Big((r(x)\wedge r(y))\int_0^1(1-z )q(\d z)-\beta_0\Big)} \cr
\ar\ar\qquad\qquad\qquad\qquad\qquad\qquad +{ \frac{1}{2}\lambda_0(x - y)E(\lambda_0,l_0)}-
 2 l_0 \frac{|r(x) - r(y)|}{x-y}\bigg].
\eeqlb
\elem

\proof
Notice that $\psi(x - y) \le 2 \le \theta \le \phi(x).$ By \eqref{L_dlel}, for $x  \le x_0$,
\beqnn
\tilde{L}f(x,y)\ar\le\ar  \psi(x - y)\Big[{(\gamma(x)+\beta_0 x) \phi'(x)+ x(\sigma^2 + x\beta_1^2)\frac{3}{x_0^2}} \cr
\ar\ar\qquad\qquad\quad +  {\frac{3x}{x_0^2}\int_0^1 z^2 \mu(\d z)  + \frac{3x}{x_0 }\int_1^\infty z  \mu(\d z)}\cr
\ar\ar\qquad\qquad\quad +{\frac{3x^2}{x_0^2}\int_{ -\infty}^{1}(\e^z - 1)^2\nu(\d z)  +  \frac{3x}{x_0}\int_1^\infty(\e^z - 1)\nu(\d z) }\cr
\ar\ar\qquad\qquad\quad +{ \frac{3x}{x_0}r(x)\int_0^1(1 - z)q(\d z)}\Big]\cr
\ar\ar+(\psi(x - y)- \phi(x))\frac12y\mu_{(x - y)}(\mbb{R}_+)\cr
 \ar\ar + \phi(x)\Big\{  {-(x - y)\e^{-\lambda_0 (x - y)}\Big[\Phi(\lambda_0)+\Big((r(x)\wedge r(y))\int_0^1(1-z )q(\d z)-\beta_0\Big)\lambda_0\Big]} \cr
\ar\ar\qquad\qquad + \frac12\sigma^2 y\psi''(x - y)+\frac{1}{2}(x - y)^2\psi''(x - y)E(\lambda_0,l_0)\Big\}\cr
\ar\ar +
{2\theta l_0|r(x) - r(y)|\psi'(x - y) q((0,1])}\cr
\ar\le\ar (\gamma(x)+\beta_0 x) \phi'(x) \psi(x - y)-\frac12y[(\theta- 2)\mu_{(x - y)}(\mbb{R}_+) + \sigma^2\theta \lambda_0^2\e^{-\lambda_0(x - y)}]\cr
\ar\ar\qquad+ 6\Big[\frac{x}{x_0^2}\sigma^2 + \frac{x^2}{x_0^2}\beta_1^2  + \frac{ x}{x_0^2}\int_0^1 z^2 \mu(\d z)+\frac{ x^2}{x_0^2}\int_{ -\infty}^{1}(\e^z - 1)^2\nu(\d z)\cr
  \ar\ar\qquad\qquad\quad +  \frac{ x}{x_0 }\int_1^\infty z  \mu(\d z) +\frac{ x}{x_0}\int_1^\infty(\e^z - 1)\nu(\d z)    +\frac{x}{x_0}r(x)\int_0^1(1 - z)q(\d z)\Big]\cr
 \ar\ar - \phi(x)\lambda_0\e^{-\lambda_0 (x - y)}(x - y)\bigg\{  \Big[\Phi(\lambda_0)/\lambda_0+\Big((r(x)\wedge r(y))\int_0^1(1-z )q(\d z)-\beta_0\Big)\Big] \cr
\ar\ar\qquad\qquad \qquad\qquad\qquad\qquad+\frac{1}{2}(x - y)\lambda_0E(\lambda_0,l_0)
-2 l_0\frac{|r(x) - r(y)|}{x-y} \bigg\}.
\eeqnn
The proof is completed here.
\qed

\subsubsection{Further estimation of the coupling generator based on $S_1$}\label{Estimation of coupling generator based on $S_1$}
We first give the definition of the set $S_1$, which provides a basis as a detailed division for further discussion on the coupling generator.
For the set $S_0$ defined in the previous subsection, we set
\beqlb
S_1 \ar=\ar S_0\cup (0,x_0]^2\cup (0,M]^2=S_0\cup (0,M]^2.
\eeqlb

Since $r$ is globally Lipschitz, there exists $k_0>0$ such that $|r(x)-r(y)|\le k_0|x-y|$. According to \eqref{e:need}, we can take $\lambda_0 > \lambda_3$ in \eqref{psi} such that
\beqlb\label{lambda0}
\ar\ar {  \frac{1}{4} x_0\lambda_0 E(\lambda_0,l_0)+ \tilde{\Phi}(\lambda_0)/\lambda_0   \ge 4k_0l_0~{\rm{and}}~ \tilde{\Phi}(\lambda_0)/\lambda_0   \ge 2k_0l_0}.
\eeqlb
In particular, $\lambda_0 > 0$ depends on $l_0$.

In this part, we give the further estimation of the coupling generator based on $S_1$ for the case of $x > y$. Let $H:=\sup_{(x,y)\in S_1} H(x,x_0)$ and take $K>0$. Now we consider the two cases separately as follows.

{\bf Case (1)~$(x,y)\notin S_1$}

In this case it holds that  $x>x_0$. Then, by Lemma \ref{x>x0} and \eqref{lambda0} as well as the fact that $\lambda_0\ge \lambda_3$, $$\tilde{L}f(x, y) \le  2{ H(x,x_0)}.$$

{\bf Case (2)~$(x,y)\in S_1$}

In this case, $(x,y)\in S_1$ implies that $0<x-y\le\max\{l_0,x_0,M\}\le l_0$.

For $x> x_0$, by \eqref{case1}, we have
\beqnn
\tilde{L}f(x, y)
\ar\le\ar  2{  H}- \frac12 \theta y[\sigma^2 \lambda_0^2\e^{-\lambda_0 (x - y)} + \mu_{(x - y)}(\mbb{R}_+)] \cr
\ar\ar -\theta(x - y)\lambda_0\e^{-\lambda_0 (x - y)} \Big[ \tilde{\Phi}({ \lambda_0})/\lambda_0 +\frac{1}{2}\lambda_0(x - y)E(\lambda_0,l_0)-  2 k_0 l_0 \Big].
\eeqnn

For $x\le x_0$, by \eqref{case2}, we have
\beqnn
\tilde{L}f(x, y)
\ar\le\ar { J(x,x_0)}-\frac12y[(\theta- 2)\mu_{(x - y)}(\mbb{R}_+) + \sigma^2 \theta \lambda_0^2\e^{-\lambda_0(x - y)}] \cr
\ar\ar  - \phi(x)\lambda_0\e^{-\lambda_0 (x - y)}(x - y)\Big[  { \tilde{\Phi}({ \lambda_0})/\lambda_0} +{ \frac{1}{2}\lambda_0(x - y)E(\lambda_0,l_0)}-  2 k_0 l_0 \Big],
\eeqnn
where
\beqnn
J(x,x_0)\ar:=\ar{ (\gamma(x)+\beta_0 x) \phi'(x) \psi(x - y)}\cr\ar\ar +  {  \frac{6 x}{x_0^2}\Big[\Big(\sigma^2 +  \int_0^1 z^2 \mu(\d z)\Big) + x \Big( \beta_1^2  + \int_{-\infty}^1(\e^z - 1)^2\nu(\d z)\Big)\Big]}\cr
\ar\ar + { \frac{6x}{x_0}\Big[\int_1^\infty z  \mu(\d z)+\int_1^\infty(\e^z - 1)\nu(\d z) + r(x)\int_0^1(1 - z)q(\d z)\Big]}
\eeqnn
 Recall that $\gamma(x) = \alpha - bx- g(x)$, $\alpha>0$ and $g$ is nondecreasing and continuous with $g(0) = 0$. We can choose constant  $r:=r(x_0) \in (0, 1/2]$  which is independent of $l_0$  such that, for any $x \in (0, r x_0],$
 \beqlb\label{case1iii}
J(x,x_0)\ar\le\ar |\beta_0 - b|6r + \frac{6g(rx_0)}{x_0}  - \frac{3\alpha}{x_0}\left(1 - r\right)^2\cr
 \ar\ar  + \frac{6 r}{x_0}\Big[\Big(\sigma^2 +  \int_0^1 z^2 \mu(\d z)\Big)\Big] +6 r^2\Big[ \Big( \beta_1^2  + \int_{-\infty}^1(\e^z - 1)^2\nu(\d z)\Big)\Big]\cr
 \ar\ar  + 6r\Big[\int_1^\infty z \mu(\d z)+\int_1^\infty(\e^z - 1)\nu(\d z) + \sup_{x\in[0,{x_0}/{2}]}r(x)\int_0^1(1 - z)q(\d z)\Big]\cr
 \ar \le \ar - \frac{6\alpha}{x_0}.
 \eeqlb

Now we consider the following more meticulous five cases.

(i) If $x>x_0$, $0<x-y\le \frac{x_0}{2}$ and $(x,y)\in S_1$, then we have $ y \ge \frac{x_0}{2} \ge x- y$ and
\beqnn
\tilde{L}f(x, y)
\ar\le\ar  2{  H}- \frac14 \theta x_0(\sigma^2 \lambda_0^2\e^{-\lambda_0l_0} + \delta),
\eeqnn
By taking
\beqlb\label{theta1}
\theta \ge \frac{4(2H + K)}{x_0(\sigma^2 \lambda_0^2\e^{-\lambda_0l_0} + \delta)},
\eeqlb
one sees that
\beqnn
\tilde{L}f(x, y) \le -K.
\eeqnn

(ii) For the case of $x>x_0$, $\frac{x_0}{2}<x-y\le l_0$ and $(x,y)\in S_1$, we then have
\beqnn
\tilde{L}f(x, y)
\ar\le\ar  2{ H}  -\theta\lambda_0\e^{-\lambda_0 l_0}\frac{x_0}{2} \Big[{  \tilde{\Phi}({ \lambda_0})/\lambda_0}+\frac{1}{4}x_0\lambda_0E(\lambda_0,l_0)- 2 k_0 l_0  \Big]\cr
\ar\le\ar 2{  H}  -\theta\lambda_0l_0\e^{-\lambda_0 l_0}x_0 k_0.
\eeqnn
By taking
\beqlb\label{theta2}
\theta \ge \frac{2H + K}{\lambda_0l_0\e^{-\lambda_0 l_0}x_0 k_0},
\eeqlb
we have
\beqnn
\tilde{L}f(x, y) \le -K.
\eeqnn

(iii)  If $x\le r x_0$
and $(x,y)\in S_1$, then it follows by \eqref{case1iii} that
\beqnn
\tilde{L}f(x, y)
\ar\le\ar - \frac{6\alpha}{x_0}.
\eeqnn

(iv) If $rx_0< x\le x_0$, $0<x-y\le \frac{rx_0}{2}$ and $(x,y)\in S_1$, we then have $y \ge \frac{rx_0}{2} \ge x - y$ and
\beqnn
\tilde{L}f(x, y)
\ar\le\ar 6R-\frac{\theta- 2}{2}y[\mu_{(x - y)}(\mbb{R}_+) + \sigma^2  \lambda_0^2\e^{-\lambda_0(x - y)}]\cr
\ar\le\ar 6R-(\theta- 2)\frac{rx_0}{4}(\delta + \sigma^2  \lambda_0^2\e^{-\lambda_0l_0}),
\eeqnn
where
\beqnn
R \ar=\ar |\beta_0 - b| + \frac{ g(x_0)}{x_0}\cr
 \ar\ar  + \frac{  1}{x_0}\Big[\Big(\sigma^2 +  \int_0^1 z^2 \mu(\d z)\Big)\Big] +   \Big[  \beta_1^2  + \int_{-\infty}^1(\e^z - 1)^2\nu(\d z) \Big]\cr
 \ar\ar  +   \Big[\int_1^\infty z \mu(\d z)+\int_1^\infty(\e^z - 1)\nu(\d z) + \sup_{x\in[0, {x_0}]}r(x)\int_0^1(1 - z)q(\d z)\Big].
\eeqnn
By taking
\beqlb\label{theta4}
\theta \ge \frac{4(6R + K)}{rx_0(\sigma^2 \lambda_0^2\e^{-\lambda_0l_0} + \delta)}+2,
\eeqlb
we have
\beqnn
\tilde{L}f(x, y) \le -K.
\eeqnn

(v) In the case of $rx_0< x\le   x_0$, $\frac{rx_0}{2}<x-y\le l_0$ and $(x,y)\in S_1$, by \eqref{lambda0}, we have
\beqnn
\tilde{L}f(x, y)
\ar\le\ar 6R   - \phi(x)\lambda_0\e^{-\lambda_0 (x - y)}(x - y)\Big[  {  \tilde{\Phi}({ \lambda_0})/\lambda_0} +{ \frac{1}{2}\lambda_0(x - y)E(\lambda_0,l_0)}-   2 k_0 l_0  \Big]\cr
\ar\le\ar 6R   - \phi(x)\lambda_0\e^{-\lambda_0 (x - y)}(x - y)\Big[  { r\tilde{\Phi}({ \lambda_0})/\lambda_0} + { \frac{1}{4}r\lambda_0 x_0E(\lambda_0,l_0)}\cr
\ar\ar\qquad\qquad\qquad\qquad\qquad\qquad\qquad +(1-r)\tilde{\Phi}({ \lambda_0})/\lambda_0 - {2 k_0 l_0 }\Big]\cr
\ar\le\ar 6R   - \phi(x)\lambda_0\e^{-\lambda_0 (x - y)}(x - y)\big[  4rk_0l_0 +2(1-r)k_0l_0- {2 k_0 l_0 }\big]\cr
\ar\le\ar 6R   - \theta r^2x_0k_0l_0 \lambda_0\e^{-\lambda_0 l_0}.
\eeqnn
By taking
\beqlb\label{theta5}
\theta \ge \frac{6R + K}{ r^2x_0k_0l_0 \lambda_0\e^{-\lambda_0 l_0}},
\eeqlb
we have
\beqnn
\tilde{L}f(x, y) \le -K.
\eeqnn

By \eqref{theta1}, \eqref{theta2}, \eqref{theta4}, \eqref{theta5}, we take
  \beqlb\label{theta}
\theta \ar:=\ar \max\Big\{4, \frac{4(2H + K)}{x_0(\sigma^2 \lambda_0^2\e^{-\lambda_0l_0} + \delta)}, \frac{2H + K}{\lambda_0l_0\e^{-\lambda_0 l_0}x_0 k_0}, \cr
\ar\ar\qquad\quad\frac{4(6R + K)}{rx_0(\sigma^2 \lambda_0^2\e^{-\lambda_0l_0} + \delta)}+2, \frac{6R + K}{ r^2x_0k_0l_0 \lambda_0\e^{-\lambda_0 l_0}}\Big\}
\eeqlb
in \eqref{phi}. Combining with all the estimates in the above cases, we then obtain the following statement.

\begin{prop}\label{P3.8}
	Suppose that Conditions $\ref{Ncond}$--$\ref{Ccond}$ are satisfied. Then  for $(x,y)\in S_1$,
	\beqnn
	\tilde{L}f(x, y) \le - K_0,
	\eeqnn where  $K_0=\min\{K,\frac{6\alpha}{x_0}\}$ is given in the beginning of Subsection $\ref{section detailed estimate}$.
\end{prop}

\subsection{Proofs of Theorem \ref{main result} and Theorem \ref{example}}

 Let $f \in \mcr{D}(\tilde{L})$ be given by \eqref{fxy} with the explicit $\psi$ and $\phi$ fixed in the previous subsections, and let
\beqlb\label{F}
F(x, y) = (V(x) + V(y))1_{\{x \neq y\}} + \varepsilon f(x, y), \qquad  (x, y) \in \mbb{R}_+^2
\eeqlb with $\varepsilon=3\lambda_2 K_0^{-1}$.
\begin{prop}
Suppose that Conditions $\ref{Ncond}$--$\ref{Ccond}$ are satisfied. Then $F \in \mcr{D}(\tilde{L})$ and there is a constant $\lambda > 0$ such that
\beqlb\label{final estimate}
\tilde{L}F(x, y) \le -\lambda F(x, y), \qquad (x, y) \in \Delta^c.
\eeqlb
\end{prop}
\proof
For any $(x, y) \in \Delta^c$, it follows by the definition of $\tilde{L}$ that
\beqnn
\tilde{L} (V(x) + V(y))1_{\{x \neq y\}} \le \tilde{L} (V(x) + V(y)) = LV(x) + LV(y).
\eeqnn

For $(x,y)\in S_1$, we have  $0<x-y\le l_0$ and
\beqlb\label{est of L short dist}
\tilde{L}F(x,y)\ar\leq\ar -\lambda_1(V(x)+V(y))+2\lambda_2-\varepsilon K_0\cr
\ar=\ar { -\lambda_1(V(x)+V(y))-\lambda_2}
\eeqlb
by the fact that $\varepsilon=3\lambda_2 K_0^{-1}$.

For $(x,y)\notin S_1$, one can see that $(x,y)\notin S_0$ and $x>M$, and that
\beqnn
V(x)\ge 12,
\eeqnn
as well as that
\beqnn
1-\frac{9\lambda_2}{K_0\lambda_1}\frac{H(x,x_0)}{V(x)}\ge
1-\frac{9\lambda_2}{K_0\lambda_1}\frac{H(x)}{V(x)}\ge\frac{1}{4}.
\eeqnn
Then
\beqlb  \label{est of L large dist}
\tilde{L}F(x,y)  \ar\leq\ar - \lambda_1(V(x)+V(y)) +2\lambda_2 +2\varepsilon H(x,x_0)\cr
\ar\leq\ar -\frac{2}{3} \lambda_1(V(x)+V(y))   +2\varepsilon H(x,x_0)\cr
\ar\leq\ar -\frac{2}{3} \lambda_1(V(x)+V(y))   +6\lambda_2 K_0^{-1} H(x,x_0)\cr
\ar\leq\ar -\frac{2}{3} \lambda_1(V(x)+V(y))\Big[1- \frac{9\lambda_2 }{K_0\lambda_1}\frac{H(x,x_0)}{V(x)+V(y)}\Big]\cr
\ar\leq\ar -\frac{2}{3} \lambda_1(V(x)+V(y))\Big[1- \frac{9\lambda_2 }{K_0\lambda_1}\frac{H(x,x_0)}{V(x)}\Big]\cr
\ar\leq\ar -\frac{1}{12} \lambda_1(V(x)+V(y))-\frac{1}{12} \lambda_1 V(x) \cr
\ar\leq\ar -\frac{1}{12} \lambda_1(V(x)+V(y))-\frac{1}{12} \lambda_1.
\eeqlb

Combining (\ref{est of L short dist}) with (\ref{est of L large dist}), there exists a constant $C_3>0$ such that
$$
\tilde{L}F(x,y)\leq-C_3(V(x)+V(y)+1).
$$
Notice that there exists a constant $C_4\geq1$ such that
$$
C_4^{-1}(V(x)+V(y)+1)\leq F(x,y)\leq C_4(V(x)+V(y)+1)
$$
for all $(x,y) \in \Delta^c$. Therefore, the result \eqref{final estimate} holds with $\lambda = (\lambda_1\wedge C_3)C_4^{-1} > 0$. \qed

\noindent  {\it Proof of Theorem $\ref{main result}$.} Let $(X_t, Y_t)_{t \ge 0}$ be the Markov coupling defined by (2.6) with $(X_0, Y_0) = (x, y)$. Recall that $\tau_m=\inf\{t\ge0: Y_t> m ~\mathrm{or} ~|X_t - Y_t| < 1/m\}$. Then, by \eqref{final estimate}, for all $t>0$,
\beqnn
\mbb{E}(\e^{\lambda(t\wedge\tau_m)}F(X_{t\wedge\tau_m}, Y_{t\wedge\tau_m}))\le F(x, y).
\eeqnn
Since the coupling process $(X_t, Y_t)_{t \ge 0}$ is non-explosive, we have $\tau_m\uparrow T$ a.s. as $m\rightarrow\infty$, where $T$ is the coupling time of the process $(X_t, Y_t)_{t \ge 0}$. Letting $m\rightarrow\infty$, by Fatou's lemma,
\beqnn
\mbb{E}(\e^{\lambda(t\wedge T)}F(X_{t\wedge T}, Y_{t\wedge T}))\le F(x, y).
\eeqnn
Since $F(x, x) = 0$ for $x \ge 0$ and $X_{T +t} = Y_{T +t}$ for $t \ge 0$, we have
 \beqnn
\mbb{E}(\e^{\lambda t}F(X_{t}, Y_{t}))\le F(x, y),\quad t>0,
\eeqnn
which clearly implies $F$ satisfies the {\it exponential contraction property:}
\beqlb\label{exponential contraction property}
\tilde{P}_tF(x,y)\leq e^{-\lambda t}F(x,y),\quad t>0,
\eeqlb
where $(\tilde{P}_t)_{t\geq0}$ is the transition semigroup of the coupling process $(X_t, Y_t)_{t \ge 0}$. It is easy to see that $\tilde{P}_t((x,y),\cdot)$ is a coupling of $P_t(x,\cdot)$ and $P_t(y,\cdot)$. Then by (\ref{exponential contraction property}) and the fact that
\beqlb
c_1F(x,y)\leq d_V(x,y)\leq c_2F(x,y),\quad  (x,y)\in\mbb{R}_+^2
\eeqlb
for some $c_2\geq c_1>0$ and $d_V(x,y)$ defined by (\ref{defin dV}), one can see
\beqlb\label{exponential contract}
W_{V}(P_t(x,\cdot),P_t(y,\cdot))\leq C_0\e^{-\lambda t}d_{V}(x,y),\quad t\geq0
\eeqlb
holds with $C_0=c_2/c_1$. By standard arguments (see, e.g., \cite[p.\ 31--32]{LLWZ23} or \cite[p.\ 601--602]{FJKR23}), (\ref{exponential ergodic}) follows.\qed

\bigskip

Finally, we can present the

\noindent {\it Proof of Theorem $\ref{example}$.} By Theorem \ref{main result}, it suffices to prove that
(\ref{condition of exam}) implies that Condition \ref{Lcond} is satisfied.
 Notice that
\beqnn
LV_{\theta}(x) \ar=\ar \theta\frac{\gamma(x)+\beta_0 x}{x + 1}V_{\theta}(x)
+\frac{\theta(\theta-1)}{2}\frac{\sigma^2x+\beta_1^2x^2}{(x+1)^2}V_{\theta}(x)\cr
\ar\ar + xV_{\theta}(x)\int_0^\infty \left[\left(1 + \frac{z}{x + 1}\right)^\theta - 1 - \theta \frac{z}{x + 1}\right]\mu(\d z)\cr
\ar\ar + V_{\theta}(x)\int_{-\infty}^\infty\left[\left(1 + \frac{x}{x+1}(e^z - 1)\right)^\theta - 1 - \theta\frac{x}{x + 1}(e^z - 1)\right]\nu(\d z)\cr
\ar\ar +{r(x)}V_{\theta}(x)\int_0^1\left[\left(1 + \frac{x}{x+1}(z - 1)\right)^\theta - 1\right]q(\d z)
\eeqnn
with $\gamma(x) = \alpha - bx- g(x).$

We have
\beqnn
 0 \ar\ge\ar \lim_{x \rightarrow \infty}x\int_0^1 \left[\left(1 + \frac{z}{x + 1}\right)^\theta - 1 - \theta \frac{z}{x + 1}\right]\mu(\d z)\cr
\ar=\ar \lim_{x \rightarrow \infty} \frac{\theta(\theta - 1)}{2}\frac{x}{(x + 1)^2}\int_0^1 (1 + \xi)^{\theta - 2}z^2\mu(\d z)\cr
\ar\ge\ar \lim_{x \rightarrow \infty} \frac{\theta(\theta - 1)}{2}\frac{x}{2(x + 1)^2}\int_0^1 z^2\mu(\d z) = 0,
\eeqnn
where $\xi \in (0, z/(x + 1)).$ Furthermore,
\beqnn
\lim_{x \rightarrow \infty}x \left[\left(1 + \frac{z}{x + 1}\right)^\theta - 1 - \theta \frac{z}{x + 1}\right] \ar=\ar \lim_{x \rightarrow \infty} \frac{\left[\left(1 + \frac{z}{x + 1}\right)^\theta - 1 - \theta \frac{z}{x + 1}\right]}{x^{-1}}\cr
\ar=\ar \lim_{x \rightarrow \infty}\frac{x^2}{(x + 1)^2}\left[\theta z\left(1 + \frac{z}{x + 1}\right)^{\theta - 1} - \theta z\right]\cr
\ar=\ar \lim_{x \rightarrow \infty} \left[\theta z\left(1 + \frac{z}{x + 1}\right)^{\theta - 1} - \theta z\right] = 0,
\eeqnn
and
\beqnn
x \left|\left(1 + \frac{z}{x + 1}\right)^\theta - 1 - \theta \frac{z}{x + 1}\right| \le 2\theta \frac{xz}{x + 1} \le 2\theta z.
\eeqnn
By dominated convergence theorem, it implies that
\beqnn
\lim_{x \rightarrow \infty}x\int_1^\infty \left[\left(1 + \frac{z}{x + 1}\right)^\theta - 1 - \theta \frac{z}{x + 1}\right]\mu(\d z) = 0.
\eeqnn

With all the conclusions above, the desired assertion follows.
\qed
	
\section*{Acknowledgements}
The research is supported by National Key R\&D Program of China (No. 2022YFA1006003), NSFC grant (No. 12201119, No. 12301167), the Education and Scientific Research Project for Young and Middle-aged Teachers in Fujian Province
of China (No. JAT231015), Guangdong Basic and Applied Basic Research Foundation (No. 2022A1515110986) and Shenzhen National Science Foundation (No. 20231128093607001).


\black

\bigskip

{\footnotesize\bf Shukai Chen:}
 {\footnotesize School of Mathematics and Statistics, 
Fujian Normal University, Fuzhou, P.R. China.}
   \texttt{skchen@fjnu.edu.cn}

{\footnotesize\bf Rongjuan Fang:}
 {\footnotesize School of Mathematics and Statistics, 
Fujian Normal University, Fuzhou, P.R. China.}
   \texttt{fangrj@fjnu.edu.cn}

{\footnotesize\bf Lina Ji:}
{\footnotesize MSU-BIT-SMBU Joint Research Center of Applied Mathematics, Shenzhen MSU-BIT University, Shenzhen, P.R. China.}
   \texttt{jiln@smbu.edu.cn}

{\footnotesize\bf Jian Wang:}
 {\footnotesize School of Mathematics and Statistics, 
Fujian Normal University, Fuzhou, P.R. China.}
   \texttt{jianwang@fjnu.edu.cn}

\end{document}